\documentclass[pdflatex,sn-mathphys-ay]{sn-jnl}

\usepackage{graphicx}%
\usepackage{multirow}%
\usepackage{amsmath,amssymb,amsfonts}%
\usepackage{amsthm}%
\usepackage{mathrsfs}%
\usepackage[title]{appendix}%
\usepackage{xcolor}%
\usepackage{textcomp}%
\usepackage{manyfoot}%
\usepackage{booktabs}%
\usepackage{url}
\usepackage{listings}%
\usepackage{multirow}
\usepackage{optidef}
\usepackage{hyperref}
\usepackage{placeins}
\usepackage{url}
\usepackage{cleveref}
\usepackage{tabularx} 
\usepackage{subcaption}

\graphicspath{ {./figures/} }
\usepackage{hyperref}
\usepackage{float}
\usepackage{verbatim} 
\usepackage{apalike}
\restylefloat{figure}
\floatstyle{plaintop} 
\restylefloat{table}

\usepackage{tikz}
\usetikzlibrary{calc,matrix, positioning, arrows, arrows.meta, shapes.geometric, fit}
\usepackage{adjustbox}
\usepackage{comment}

\usepackage[norelsize, ruled, vlined, boxed, commentsnumbered]{algorithm2e}

\usepackage{lineno}
\usepackage{algpseudocode}  

\begin{document}

\title[RKO for MIPs]{{\center Applying a Random-Key Optimizer\\ to Mixed Integer Programs}}



\author*[1]{\fnm{Antonio A.} \sur{Chaves}}\email{antonio.chaves@unifesp.br}

\author[1,3]{\fnm{Mauricio G.C.} \sur{Resende}}\email{{mgcresende@unifesp.br}}

\author[2]{\fnm{Carise E. } \sur{Schmidt}}\email{carise.schmidt@ifsc.edu.br}

\author[4]{\fnm{J. Kyle} \sur{Brubaker}}\email{j.kylebrubaker@gmail.com}

\author[5]{\fnm{Helmut G.} \sur{Katzgraber}}\email{hgk@55n.vc}

\author[6,7]{\fnm{Martin J.A.} \sur{Schuetz}}\email{maschuet@amazon.com}

\affil*[1]{\orgname{Federal Univ. of São Paulo}, \orgaddress{\city{S. J. dos Campos},  \state{SP}, \country{Brazil}}}

\affil[2]{\orgname{Instituto Federal de Santa Catarina}, \orgaddress{\city{Chapec\'o},  \state{SC}, \country{Brazil}}}

\affil[3]{\orgname{{DIMACS}}, \orgaddress{\city{{Piscataway}},  \state{{NJ}}, \country{USA}}}

\affil[4]{\orgaddress{\city{{San Diego}},  \state{{CA}}, \country{USA}}}

\affil[5]{\orgname{{55 North}}, \orgaddress{\city{{Copenhagen}}, \country{Denmark}}}

\affil[6]{\orgdiv{Amazon Advanced Solutions Lab}, \orgaddress{\city{Seattle},  \state{WA}, \country{USA}}}

\affil[7]{\orgdiv{AWS Center for Quantum Computing}, \orgaddress{\city{Pasadena},  \state{CA}, \country{USA}}}


\abstract{
Mixed-Integer Programs (MIPs) are NP-hard optimization models that arise in a broad range of decision-making applications, including finance, logistics, energy systems, and network design. Although modern commercial solvers have achieved remarkable progress and perform effectively on many small- and medium-sized instances, their performance often degrades for large-scale or highly constrained formulations. This paper explores the use of the Random-Key Optimizer (RKO) framework as a flexible, metaheuristic alternative for computing high-quality solutions to MIPs through the design of problem-specific decoders.
The proposed approach separates the search process from feasibility enforcement by operating in a continuous random-key space while mapping candidate solutions to feasible integer solutions via efficient decoding procedures. We evaluate the methodology on two representative and distinct benchmark problems: the mean–variance Markowitz portfolio optimization problem with buy-in and cardinality constraints, and the Time-Dependent Traveling Salesman Problem. For each formulation, tailored decoders are developed to reduce the effective search space, promote feasibility, and accelerate convergence. Computational experiments demonstrate that RKO consistently produces competitive, and in several cases superior, solutions compared to a state-of-the-art commercial MIP solver, both in terms of solution quality and computational time. These results highlight the potential of RKO as a scalable and versatile heuristic framework for tackling challenging large-scale MIPs.

}




\maketitle

\setcounter{algocf}{0}  

\vspace{1.0cm}
\section{Introduction}
\label{introduction}

Mixed Integer Programs (MIPs) encompass NP-hard optimization problems where the goal is to minimize some objective subject to constraints, with some or all of the variables constrained to be integer-valued \citep{Karp1972, papadimitriou1998combinatorial}. Applications of MIPs can be found in virtually every industry, in areas such as transportation and logistics, telecommunications, manufacturing, and finance. Prominent examples include capacity planning, resource allocation, bin packing, and portfolio optimization, among others \citep{korte2011combinatorial}. Commercial black-box solvers such as CPLEX \citep{IBMCPLEX2024}, Gurobi \citep{gurobi}, or Xpress \citep{FICOXpress2025} are typically used to solve MIPs, leaving end users with potentially costly licensing fees and limited insights into the anatomy of the underlying algorithms used. 
Moreover, these solvers typically rely on the branch-and-bound (B\&B) paradigm that can solve small and intermediate-size problems efficiently to optimality, but often suffer from exponential runtimes in the size of the input, calling for the development of novel, efficient methods that can deal with the ever-increasing scale of real-world optimization problems. 

\begin{figure}[ht]
    \centering
    \includegraphics[width=0.9\linewidth]{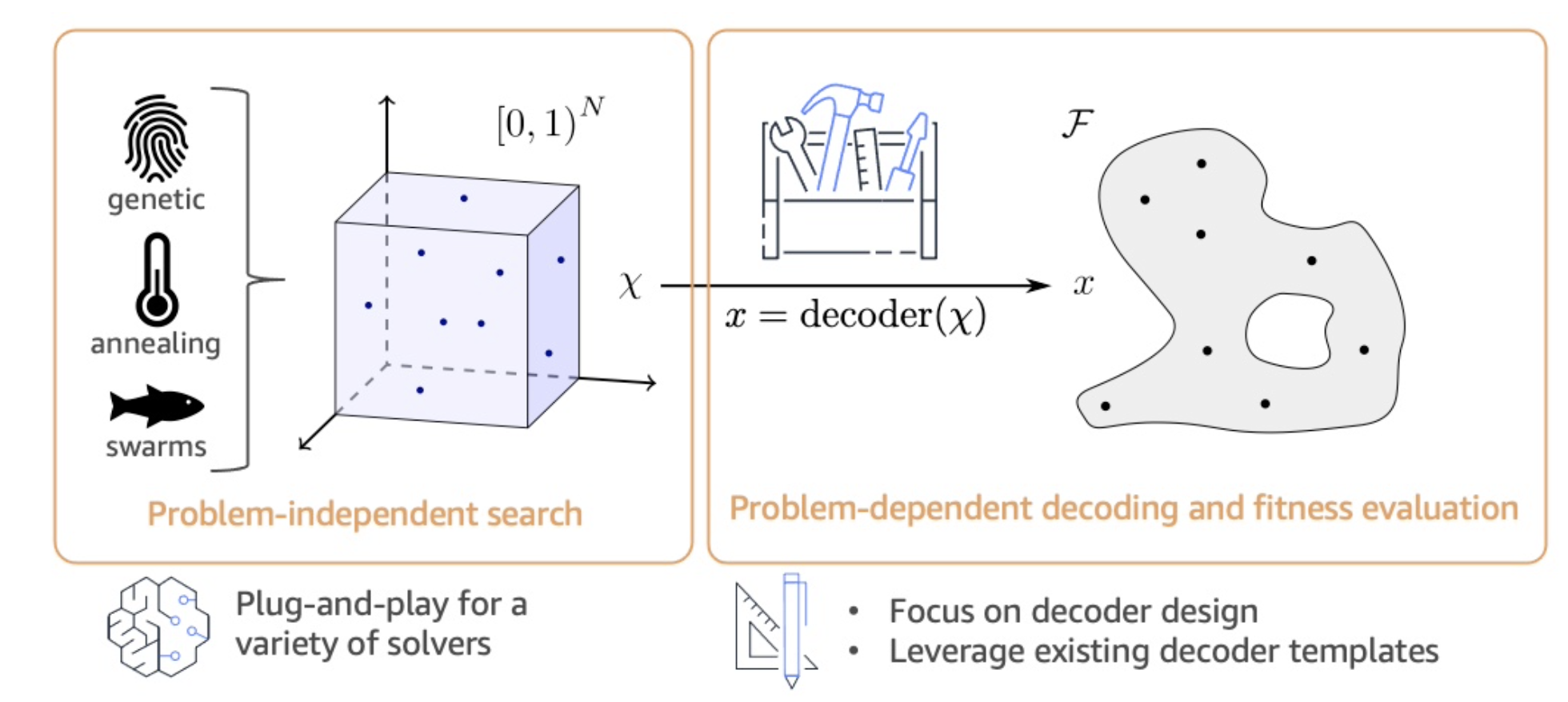}
    \caption{Schematic illustration of our approach. Our framework is based on a distinct separation of problem-independent (blue hypercube) and problem-dependent (decoder) modules. The search for high-quality solutions is performed within the (problem-independent) space of random keys, as can be implemented with a stack of algorithmic paradigms such as genetic algorithms, simulated annealing, or swarm-based algorithms. The deterministic (problem-dependent) decoder maps a given random-key vector $\mathcal{X} \in [0,1)^n$ to a feasible solution $x=\mathit{decoder(\mathcal{X})}$ within the problem space (as displayed by black dots). For mixed integer programs, our decoder design allows for an efficient, native encoding of a plethora of common constraints such as integrality constraints, variable bounds, and cardinality constraints.}
    \label{fig:SchematicRKO}
\end{figure}

This work shows how the Random-Key Optimizer (RKO) \citep{Chaves2025-RHO-JH} can be used to solve MIPs. Specifically, we outline a general method that can leverage a plethora of optimization paradigms such as genetic algorithms, physics-inspired algorithms (e.g., simulated annealing), or swarm-based algorithms (e.g., particle swarm optimization, or ant colony optimization), all within one unified framework, while natively encoding a large class of constraints through problem-specific decoder design, such as integrality constraints, variable bounds, and cardinality constraints of the MIPs. Our approach is schematically depicted in \autoref{fig:SchematicRKO}. 
Next, we illustrate our approach with numerical experiments for two benchmark problems: the canonical mean-variance Markowitz portfolio optimization problem in the presence of buy-in and cardinality constraints \citep{chang2000heuristics} and the Time-Dependent Traveling Salesman Problem (TD-TSP) \citep{malandraki1992time}.

The key contributions of this study include:

\begin{itemize}
    \item We developed a continuous optimization algorithm designed for the efficient solution of MIPs;

    \item We created innovative decoders that directly map random-key solutions to the Markowitz portfolio and the TD-TSP models;

    \item Our RKO method demonstrates high-quality solutions with competitive or superior computational performance compared to a state-of-the-art commercial solver.
\end{itemize}

The structure of this paper is as follows. Section \ref{sec:sec2} introduces the RKO. Section \ref{sec:sec3} establishes the theoretical framework for MIPs. Section \ref{sec:sec4} details the encoding/decoding scheme for both the Mixed-Integer Limited-Asset Markowitz and the TD-TSP models. Section \ref{sec:sec5} presents our experiments and results. Finally, Section \ref{sec:sec6} concludes with a summary of the main contributions and discusses the implications of this research.


\section{Random-key Optimizer - RKO}
\label{sec:sec2}

Our approach utilizes the concept of random keys, as proposed in the Biased Random-Key Genetic Algorithm (BRKGA) \citep{gonccalves2011biased}. To establish our notation and terminology, we begin with a brief review of BRKGA and demonstrate how this algorithm can be generalized to alternative solution strategies. A refinement of the RKGA proposed by \citet{Bean_RKGA_1994}, BRKGA features a heuristic structure for solving optimization problems. Although most of the work in the literature has focused on combinatorial optimization problems, BRKGA has also been applied to continuous optimization problems \citep{silva2014finding}. BRKGA is based on the idea that any solution to an optimization problem can be encoded as a vector of random keys, i.e., a vector $\mathcal{X}$ where each entry is a real number, randomly generated in the interval $[0, 1)$. Such a vector $\mathcal{X}$, inherently defined in the hypercube $[0,1)^n$, is mapped to a feasible solution for the combinatorial optimization problem using the decoder. This deterministic algorithm takes a vector of random keys as input and returns a feasible solution to the optimization problem and the solution's cost. A penalty is added to the cost in case of infeasibility.

The BRKGA has been generalized into a general RKO paradigm \citep{Chaves2025-RHO-JH}. Similar to BRKGA, RKO uses a vector of random keys to encode a solution for a combinatorial optimization problem and employs a decoder to evaluate solutions encoded by the random-key vector. However, in RKO, the problem-independent search in the random-key space is generalized beyond the genetic algorithm towards a much broader class of metaheuristics, using, e.g., Genetic Algorithm (GA), Simulated Annealing (SA), Greedy Randomized Adaptive Search Procedure (GRASP), Iterated Local Search (ILS), Variable Neighborhood Search (VNS), Particle Swarm Optimization (PSO), and Large Neighborhood Search (LNS).

The RKO framework has been extensively expanded to integrate various metaheuristic paradigms. For instance, PSO algorithms have been adapted to work with random-key encoding for permutation generation, as presented by \citet{lin2010efficient} and \citet{bewoor2017evolutionary,bewoor2018production}. Similarly, \citet{garcia2015random} devised a Harmony Search (HS) heuristic that leverages random-key encoding to effectively tackle job-shop scheduling problems, ensuring the generation of feasible scheduling solutions through the application of HS operators to random-key encoded harmonies. \citet{ouaarab2015random} present the random-key Cuckoo Search algorithm to solve the traveling salesman problem, in which new solutions are generated via Lévy flights \citep{yang2009cuckoo}.


Further contributions include the work of \citet{pessoa2018heuristics}, who introduced four constructive heuristics for the flowshop scheduling problem, alongside BRKGA, ILS, and Iterated Greedy Search (IGS) approaches that explore sequence-based neighborhoods. Similarly, \citet{andrade2019scheduling} developed BRKGA and ILS approaches for machine-dependent scheduling problems, employing a decoder that converts random-key vectors into schedules. This adaptable permutation-based decoder was utilized for Tabu Search (TS) and SA implementations.

Innovations like implicit path-relinking (IPR), introduced by \citet{andrade2021multi}, further extend the RKO's capabilities by exploring the separation between problem and solution spaces through path construction within the unit hypercube, using the decoder for solution evaluation. The RKO framework also found application in robot motion planning, with \citet{Schuetz_RKO_2022} combining BRKGA with dual annealing for efficient random-key space exploration. This work also formally introduced the term ``RKO'', resulting in a patent for the authors \citep{SchBruResKat2025a}.

More recently, RKO implementations based on SA, ILS, and VNS have been deployed for the tree hub location problem by \citet{ManPolMacJulProGiaSalCha23a}. Expanding the metaheuristic portfolio, \citet{RK-GRASP} developed a GRASP-based RKO approach, showcasing its effectiveness across a diverse range of optimization challenges, including the traveling salesman problem, Steiner triple covering problem, node capacitated graph partitioning problem, and the job sequencing problem with tool switching.

The general pseudo-code for RKO is described in Algorithm \ref{alg:rko}. It starts by accepting as input the instance that needs to be solved and a computation time limit ($\mathit{run}_{\mathit{time}}$). Alternatively, one could use as a stopping rule a maximum number of calls to the decoder ($\mathit{run}_{\mathit{dcalls}}$). Next, a user-specified set of metaheuristic algorithms is established. Lines 2 and 3 read the instance data and initialize the pool of solutions at random. Each metaheuristic algorithm runs concurrently in the parallel procedure described in lines~4-6, exchanging information through the pool of solutions. Finally, the algorithm returns the best solution from this pool.

\begin{algorithm}[htbp]
\small
\DontPrintSemicolon
\KwIn{\textit{instance}, $\mathit{run_{time}}$}
\KwOut{$\mathcal{X}_{\mathit{best}}$}
\Call{MH}{} $\leftarrow \{\mathit{MH}_1, \mathit{MH}_2, \dots, \mathit{MH}_m\}$ \tcp*[f]{Create metaheuristic list} \;
\Call{ReadData}{$\mathit{instance}$} \tcp*[f]{Read instance data} \;
\Call{CreatePoolSolutions}{$\mathit{pool}[ \ ], \mathit{size_{pool}}$} \tcp*[f]{Random solution pool} \;
\While{$\mathit{currentTime} \le \mathit{run_{time}}$}{
    \For(\tcp*[f]{Run all MHs in parallel}){$i \in$ MH \textbf{in parallel}} {                     

        \Call{MH}{$i$}

    }
}
$\mathcal{X}_{\mathit{best}} \leftarrow \mathit{pool}[0]$
\caption{Random-Key Optimizer (RKO)}
\label{alg:rko}
\end{algorithm}

Several key components form the foundation of the RKO. It has a pool of elite solutions that carefully preserves the best and most varied solutions found during the optimization process, as well as a generator of random initial solutions to start the search process. Four different moves are employed by a shaking mechanism to perturb random-key vectors, enabling efficient exploration of the solution space. Furthermore, two solutions are combined to create a new one using a blending method, which is an adaptation of uniform crossover. Lastly, to direct the search for the best solutions, the framework uses a variety of metaheuristic algorithms in addition to local search heuristics.

This method starts with $size_{pool}$ solutions and maintains a shared pool of elite solutions. Any new solution found by a metaheuristic during the search process that outperforms its current best is considered for inclusion in this elite pool. This inclusion is subject to two restrictions: the new solution cannot be a copy of any existing solution in the pool, and if it is approved, it replaces the solution with a worse objective value that is most similar to it, ensuring that the pool's size remains constant.

To add controlled perturbations to a random-key vector, the shaking method, which is essential for diversification, first chooses a perturbation rate $\beta$ at random.  Four different neighborhood moves are then applied based on this rate. Two randomly selected keys, $\mathcal{X}[i]$ and $\mathcal{X}[j]$, have their values switched by the Swap move. A key $\mathcal{X}[i]$ is swapped with its immediate successor, $\mathcal{X}[i+1]$, using the Swap Neighbor move. By using the Mirror move, a key $\mathcal{X}[i]$ is replaced by its complement $1.0-\mathcal{X}[i]$. Moreover, the Random move ensures a thorough exploration of the solution space by substituting a completely new random value for a key $\mathcal{X}[i]$.

The blending method combines two solutions, $\chi_a$ and $\chi_b$, to create a new random-key vector $\chi_c$. By incorporating additional stochastic components, this procedure builds upon the concept of Uniform Crossover (UX) \citep{davis1989handbook}. There is a probability $\rho$ of $\chi_c$ directly inheriting a value from either solution $\chi_a$ or solution $\chi_b$ for every key $\chi_c[i]$. The contribution from solution $\chi_b$ is further adjusted by a factor parameter: a factor of $1$ uses $\chi_b[i]$, while a factor of $-1$ uses $1.0 - \chi_b[i]$. Additionally, to increase the diversity of the generated solutions, a new random value is assigned to a key with a probability $\mu$.

To exploit the search process, we utilize Randomized Variable Neighborhood Descent (RVND) \citep{subramanian2010parallel}, which dynamically randomizes the order of applying local search heuristics in each iteration.  To achieve this, we developed four problem-independent local search heuristics specifically designed to operate in the random-key solution space: Swap LS, Mirror LS, Farey LS, and an Adapted Nelder-Mead LS \citep{nelder1965simplex}.  Each of these heuristics evaluates the quality of neighbor solutions using a decoder.

Beyond traditional evolutionary algorithms, the RKO framework supports a wide range of metaheuristics designed especially for the random-key solution space. Before using its search method, each metaheuristic
algorithm in the RKO begins its search with randomly generated solutions. Nine classical metaheuristics are currently included in the RKO framework. In addition to an Implicit Path-Relinking (IPR) mechanism, there are the following: SA \citep{Kirkpatrick_SA_1983}, GRASP \citep{Feo_GRASP_1995}, ILS \citep{Lourenco_ILS_2003}, VNS \citep{Mladenovic_VNS_1997}, LNS \citep{pisinger2010}, PSO  \citep{Kennedy_PSO_1995}, GA \citep{holland1992adaptation}, BRKGA \citep{gonccalves2011biased}, and BRKGA with Clustering Search (BRKGA-CS) \citep{Chaves_BRKGA_QL_2021}.

\cite{Chaves2025-RHO-JH} provide a detailed explanation of each of these components as well as the operations of the RKO framework. The RKO framework is publicly available in the GitHub repository\footnote{\url{https://github.com/RKO-solver}}, including its C++ source code and extensive documentation, for more accessibility and transparency.


\section{General Mixed Integer Programs}
\label{sec:sec3}

In this section, we focus on the class of linear mixed-integer programs (MILPs). Applications of MILPs include production planning (e.g., job-shop scheduling), scheduling (e.g., vehicle scheduling in transportation networks), and telecommunication network design, among others.

Generically, MILPs can be formalized in compact form as
\begin{equation}
\underset{x}{\operatorname{argmin}} \{c^Tx \mid Ax \leq b, l \leq x \leq u, x \in \mathbb{Z}^p \times \mathbb{R}^{n-p} \}
\label{eq:milp}
\end{equation}
where the cost vector $c \in \mathbb{R}^n$ specifies the objective function, the matrix $A \in \mathbb{R}^{m \times n}$ together with the vector $b \in \mathbb{R}^m$ describes the constraints, the vectors $l, u \in \mathbb{R}^n$ give the lower and upper variable bounds, respectively, and $\mathcal{I} = \{1, ..., p\} \subseteq  [n]$ refers to the index set of integer variables. For $p = n$, we recover an integer linear program (ILP) known to be NP-hard \citep{Karp1972}. The latter reduces to 0/1 programming if all variables are binary. The size of the problem is typically quantified by the number of variables $n$ and the number of constraints $m$. While any point $x \in \mathbb{R}^n$ is a (complete) assignment, feasible assignments are those that satisfy all the constraints specified in Eq. \eqref{eq:milp}, and an optimal assignment (or solution) refers to a feasible assignment that also minimizes the objective \citep{boyd2004convex}.

\subsection{Theoretical Framework}

We now show how any random-key algorithm can be applied to mixed-integer programs, including MILPs. To this end, we propose a decoder design that natively accounts for integrality constraints, $x \in \mathbb{Z}^p \times \mathbb{R}^{n-p}$, as well as variable bounds, $l \le x \le u$. Specifically, for any random key $\mathcal{X}_i \in [0,1)$ we simply take

\begin{equation*}
x_i=
\begin{cases}
\text{int}[l_i - 0.5 + (u_i - l_i + 1) \mathcal{X}_i], & \text{if } i \in \mathcal{I}; \\
l_i + (u_i - l_i) \mathcal{X}_i, & \text{otherwise,}
\end{cases}
\label{eq:xi}
\end{equation*}

\noindent with $int[\cdot]$ rounding to the closest integer. This mapping yields an integer (bounded) decision variable $x_i \in \{l_i, l_i + 1,..., u_i\}$ for $i \in \mathcal{I}$, and a continuous (bounded) decision variable $x_i \in [l_i,u_i]$, as desired. For violations of the remaining constraints, $Ax \le b$, we introduce prohibitively large penalty values, allowing the algorithm to iteratively steer away from these low-fitness solutions over the course of the evolution. Because we are not bound to differentiable loss functions, we can handle these cases with a simple if-else logic. Specifically, we set $\Delta = b - Ax$, with $\Delta_i \ge 0$ indicating satisfied constraints and $\Delta_i < 0$ signaling violations thereof. For any $\Delta_i < 0$ we then add some penalty $\phi(\Delta_i)$ to the cost (fitness) value $f(x)$ as

\begin{equation}
f(x) = c^T x + \sum_{i=1}^{m} \phi(\Delta_i) \mathbf{1}_{\Delta_i < 0}
\label{eq:fx}
\end{equation}

\noindent with indicator function $\mathbf{1}_x$. The cost value $f(x)$ provides feedback to the optimizer, and is used to steer the optimizer towards feasible assignments over the course of the search process. As a simple model for the penalty function $\phi(\cdot)$ we can take $\phi(\Delta_i) = P \cdot \Delta^2_i$ with pre-factor $P > 0$, but other models such as logistic barrier functions (as used in the context of interior-point methods \citep{boyd2004convex}) may be used as well. The numerical value for the penalty parameter $P$ can be optimized in an outer loop. Overall, the proposed decoder design is summarized in pseudo-code in Algorithm \ref{alg:decoderExample}. In Section \ref{sec:sec4}, we discuss additional decoder designs that address common cardinality constraints, sub-route elimination, integrality constraints, and variable bounds.

\begin{algorithm}[htbp]
\small
\DontPrintSemicolon
\KwIn{$c, A, b, l, u, p, \mathcal{X}$}
\KwOut{$x, \mathit{cost}$}
\tcp{Initialize assignment}
$x \leftarrow \text{np.zeros}(\text{length}(\mathcal{X}))$ \;
\tcp{Get integer assignments}
\For{$i \leftarrow 1$ \KwTo $p$}{
    $x[i] \leftarrow \text{int}(l[i] - 0.5 + (u[i] - l[i] + 1) \cdot \mathcal{X}[i])$ \;
}
\tcp{Get continuous assignments}
\For{$i \leftarrow p+1$ \KwTo $\text{length}(x)$}{
    $x[i] \leftarrow l[i] + (u[i] - l[i]) \cdot \mathcal{X}[i]$ \;
}
\tcp{Initialize cost of assignment $f(x)$}
$cost \leftarrow c^T x$ \;
\tcp{Add penalties for constraint violations}
$\Delta \leftarrow b - A x$ \;
\For{$\Delta_i \in \Delta$}{
    \If{$\Delta_i < 0$}{
        $\mathit{cost} \leftarrow \mathit{cost} + \phi(\Delta_i)$ \;
    }
}
\caption{Example decoder for generic MIPs}
\label{alg:decoderExample}
\end{algorithm}


\section{Applications}
\label{sec:sec4}

In this paper, we present two case studies of RKO applications for solving MIPs, namely the Mixed-Integer Limited-Asset Markowitz Model and the Mixed-Integer TD-TSP Model. Representing the random-keys vector for these models can include all decision variables. However, this will increase execution time and unnecessarily restrict convergence. Therefore, in these approaches, we develop decoders that find values only for non-zero decision variables.

\subsection{Mixed-Integer Limited-Asset Markowitz Model}
\label{sec:sec4a}

First, to illustrate our approach, we consider the canonical Markowitz (mean-variance) approach to portfolio selection in the presence of cardinality constraints that limit a portfolio to have a specific number of assets and bounds on the proportion of the portfolio held in a given asset (if any of the asset is held) \citep{chang2000heuristics, cesarone2011portfolio} as to limit risk exposure and transaction costs for small investments through upper and lower (buy-in) constraints, respectively. This example shows that our approach is not limited to MILPs but can handle non-linear MIPs within the same framework.

Let $n$ be the total number of assets available, $\mu_i$ be the expected return of asset 
$i \in \{1,...,n\}$, $\Sigma_{ij}$ be the covariance between assets $i \in \{1,...,n\}$ and $j \in \{1,...,n\}$, and $\lambda$ be the risk-aversion parameter, which controls how much an investor penalizes variance relative to expected return, where $0 \le \lambda \le 1$. Furthermore, let $K$ be the desired number of assets in the portfolio, $l_i$ be the minimum proportion that must be held of asset $i = \{1,...,n\}$ (if any of asset $i$ is held), $u_i$ be the maximum proportion that can be held of asset $i=\{1,...,n\}$ (if any of asset $i$ is held), with $0 \le l_i \le u_i \le 1$. In practice, the lower bounds $l_i$ represent a min-buy or minimum transaction level for asset $i$, while the upper bounds $u_i$ limit the exposure of the portfolio to individual assets. The decision variables are the allocation weights $w_i$ for asset $i \in \{1,...,n\}$, with $0 \le w_i \le 1$, together with the binary indicator variables $z_i \in \{0, 1\}$, with $z_i = 1$, if asset $i$ is held, and $z_i = 0$ otherwise. We then consider the following mixed-integer optimization problem:

\begin{equation}
\min_{\mathbf{w}} \ H(\mathbf{w}) = \lambda \mathbf{w}^\top \Sigma \mathbf{w} - (1 - \lambda) \boldsymbol{\mu}^\top \mathbf{w} \label{fo}
\end{equation}

\vspace{-0.3cm}

\begin{align}
\text{s.a.} \quad & \sum_{i=1}^{n} w_i = 1, \label{r1}\\ 
& \sum_{i=1}^{n} z_i = K, \label{r2} \\ 
& l_i z_i \leq w_i \leq u_i z_i, \quad i = 1, \ldots, n, \label{r3} \\
& z_i \in \{0, 1\}, \quad i = 1, \ldots, n.  \label{r4}
\end{align}

In the limiting case where $\lambda = 0$, the objective function \eqref{fo} aims to maximize expected return (irrespective of the risk involved), and the optimal solution will involve a single asset with the highest return. Conversely, for $\lambda = 1$ our goal is to minimize risk (irrespective of the return involved), and the optimal solution will typically involve a larger number of assets. Intermediate values with $0 < \lambda < 1$ describe an explicit trade-off between risk and return. For example, for $\lambda = 0.25$ the objective becomes $0.25[risk] - 0.75[return]$. This problem features four (families of) constraints. Constraint \eqref{r1} represents the budget constraint, given that the available funds are limited. Constraint \eqref{r2} refers to the cardinality constraint, forcing us to select the best portfolio with investments in exactly $K$ assets. Note that we have explicitly chosen to formulate this problem with an equality (rather than an inequality) with respect to the number of assets $K$ in the portfolio. This is because if we can solve the equality-constrained case (as considered here), then any situation involving inequalities (lower or upper limits on the number of assets in the portfolio) can be solved in a series of equivalent problems (that are trivially parallelized). Specifically, the generalized constraint $K_L \le \Sigma_i z_i \le K_U$ can be dealt with by simply studying all valid values of $K$, i.e., $K = \{K_L,K_L + 1,...,K_U\}$ \citep{chang2000heuristics}. Constraints \eqref{r3} describe the variable bounds, requiring that all non-zero investments lie within the investment band described by lower $l$ and upper bounds $u$, respectively. These constraints couple the integer variables $z_i$ with the continuous variables $w_i$. Finally, constraints \eqref{r4} impose integrality on the variables $z_i$, making them binary indicator variables.

\subsubsection{Encoding/Decoding}

We now illustrate the simplicity and flexibility of our approach, with one decoder design for the non-linear MIP described in Eq. \eqref{fo}--\eqref{r4}. This decoder handles integrality constraints ($z_i \in \{0,1\}$), budget, and cardinality constraints natively by design, meaning that these constraints will always be fulfilled automatically by the decoder. The remaining variable bounds constraints are accounted for via soft penalty terms that effectively fold violations of these constraints into the cost function, thereby steering the solver towards feasible (low-cost) solutions throughout the search process. 


Each solution corresponds to a random-key vector of size $2 \times K$, with each random key representing a real value in the interval $[0,1)$. The decoder maps each solution into a set of selected assets and continuous weights. The first $K$ random keys are converted into an identifier corresponding to one of the available assets (an integer number belonging to the set $\{1,..., n\}$), while the second block represents the continuous weights. 

We first create a list of size $n$ containing all available assets enumerated by an identifier $1,...,n$. This decoder is applied iteratively and without replacement by selecting an asset in each of $K$ iterations. 
Initially, the first random key is multiplied by $n$ and the ceiling function is applied. This integer number obtained corresponds to the first selected asset's identifier ($id$), giving $z_{id} = 1$. This asset will be removed from the list for subsequent iterations. 
Thus, the set of assets to be considered in the next iteration has a size $n-1$.
Next, the second selected asset is selected from the second random key, multiplied by $n-1$, and the ceiling function is applied. This number corresponds to the identifier of the second selected asset. This process is repeated until $K$ assets are selected. 
After selecting the $i$-th asset, the continuous weight of the asset $id$ is set by a simple affine transformation using the random key $i+K$,  $w_{id} = l_i + (u_i - l_i) \times  \mathcal{X}[i+K]$.

\begin{figure}[htbp]
    \centering
    \includegraphics[width=1.0\linewidth]{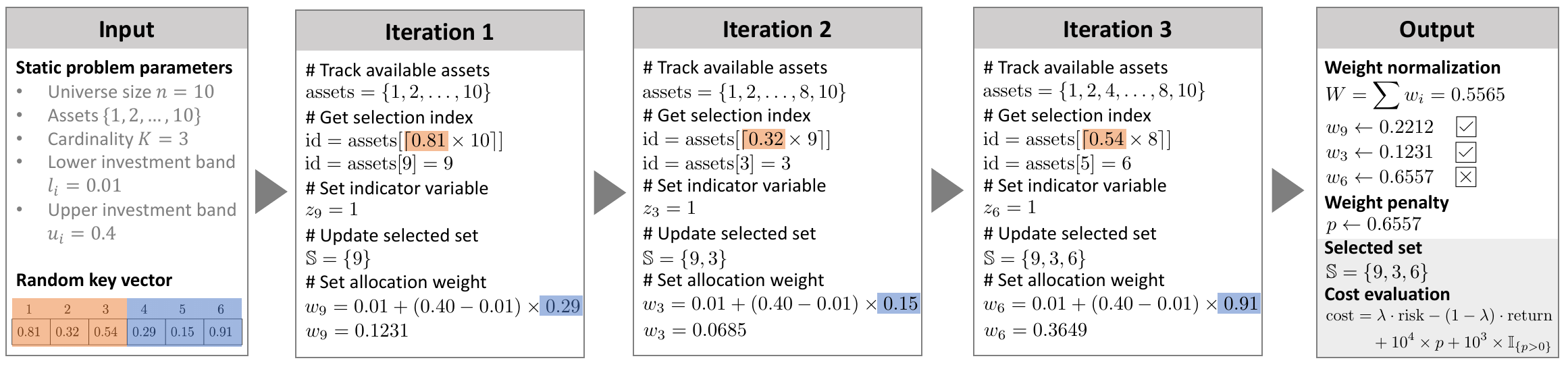}
    \caption{Decoder example for the mixed-integer limited-asset Markowitz model. The solution is encoded in the random key vector of size $2K$, with the first block encoding the integer variables $z_{i}$ and the second block representing the continous weights $w_{i}$. Decoding is done through $K$ simple iterations. For the given sample random key vector we obtain a set of selected assets given by $\mathbb{S}=\{9,3,6\}$. Cardinality, budget, and integrality constraints are satisfied by virtue of this decoder design, while the remaining variable bounds are handeled via appropriate penalty terms in the cost function.}
    \label{fig:exampleDecoder}
\end{figure}

Finally, a simple post-processing step always satisfies the budget constraint with the transformation $w_i = w_i/ \sum_{i=1}^{n} w_i$. However, this transformation may violate the variable bound constraints, which can be handled with a simple penalty term proportional to the extent to which each variable $w_i$ violates its lower and upper bounds. Fixed and variable penalties are also added to avoid overly small penalties. Figure \ref{fig:exampleDecoder} presents an example of mapping the random-key vector $\mathcal{X} = \{0.81, 0.32,$ $ 0.54, 0.25,$ $0.15, 0.91\}$ of size $K = 3$ for a dataset with $n = 10$ assets, and Algorithm \ref{alg:decoder} illustrates the pseudocode of the decoder.

\begin{algorithm}[htbp]
\small
\DontPrintSemicolon
\KwIn{$\mu, \Sigma, \lambda, l, u, K, \mathcal{X}$}
\KwOut{$w, z, cost$}
$\mathit{assets} \leftarrow \{1,...,n\}$ \;
\For{$i \leftarrow 1$ \KwTo $K$}{
    $\mathit{id} \leftarrow \mathit{assets}[\text{ceiling}(\mathcal{X}[i] \times (n+1-i)]$ \;
    $z[id] \leftarrow 1$\; 
    $w[id] \leftarrow l[\mathit{id}] + (u[\mathit{id}] - l[\mathit{id}]) \times  \mathcal{X}[i+K]$\; 
    $\mathit{assets} \leftarrow \mathit{assets} \backslash \mathit{id}$
}
\For{$i \leftarrow 1$ \KwTo $n$}{
    $w[i] \leftarrow w[i] / \sum_{i}{w_i}$\; 
    $\mathit{\mathit{penalty}} \leftarrow \max(0, w[i] - u_i) + \max(0, l_i - w[i])$ \;
}
$cost \leftarrow \lambda \mathbf{w}^\top \Sigma \mathbf{w} - (1 - \lambda) \boldsymbol{\mu}^\top \mathbf{w} + (10^4 \times \mathit{penalty} + 10^3 \times 1_{\{\mathit{penalty} > 0\}})$
\caption{Decoder for portfolio optimization problem}
\label{alg:decoder}
\end{algorithm}





\subsection{Mixed-Integer TD-TSP Model}
\label{sec:sec4b}

As a second example, we focus on the TD-TSP, in which variations in travel times over the planning horizon are explicitly taken into account in the routing decisions. A single vehicle leaves a depot, serves each customer exactly once, and returns to the depot, where the time needed to travel from one location to the next depends on the departure time and is modeled by discretizing the planning horizon into time intervals. This setting captures variations in travel conditions, such as recurrent congestion patterns, and leads to a time-indexed MILP formulation in which routing and timing decisions are tightly integrated.

Let $G=(V,A)$ be a directed graph, where $V=\{0,1,2,...,n, n+1\}$ denotes the set of nodes and $A=\{(i,j) \in V \times V: i \ne j\} \setminus \{(0,n+1),(n+1,0)\}$ denotes the set of arcs. Node $0$ represents the depot, while $n+1$ is a dummy node, i.e., a copy of the depot used to mark the end of the route. Let $V_c$ denote the subset of customers that must be served, i.e., $V_c = V \setminus \{0, n+1\}$. Each customer $i \in V_c$ has an associated deterministic service time $s_i>0$; by convention, we set $s_0=s_{n+1}=0$.

The graph $G$ is time-dependent, meaning that as traffic conditions change, the travel time on each arc $(i,j) \in A$ can vary. The planning horizon is discretized into a set $H=\{0,1,2,...,m\}$ of $m+1$ time intervals of equal length $\bar{T}$. For each time interval $h \in H$, traffic conditions are assumed to be constant over $[h \bar{T}, (h+1)\bar{T})$. Accordingly, for every arc $(i,j) \in A$ and time interval $h \in H$ let $t_{ij}^h$ denote the deterministic travel time required to traverse arc $(i,j) \in A$ when departure takes place in interval $h \in H$.

We use a formulation with a time-discretized representation of travel times \citep{malandraki1992time}, embedded in a two-commodity flow model on an extended network \citep{Baldacci2003, Baldacci2004}, and define the following decision variables. Let $x_{ij}^h$ denote binary variables that take value $1$ if arc $(i,j) \in A$ is traversed in time interval $h \in H$, and $0$ otherwise. Let $y_{ij}$ denote continuous nonnegative flow variables that support a two-commodity flow representation of the route: whenever arc $(i,j) \in A$ belongs to the route, $y_{ij}$ represents the flow associated with customers that still need to be served, whereas  $y_{ji}$ represents the flow associated with customers that have already been served. Finally, let $a_i$ denote continuous nonnegative variables representing the accumulated time upon departure from node $i \in V \setminus \{n+1\}$ or upon arrival at the terminal node $n+1$.
We consider the following mixed-integer optimization problem:
\begin{equation}
	Min \sum_{(i,j) \in A} \sum_{h \in H} x_{ij}^h t_{ij}^h
	\label{R0}
\end{equation}

\begin{align}
\text{s.a.} \quad & x_{i0}^h = 0, \qquad \forall i \in V_c, \quad \forall h \in H, \label{R1}\\ 
& x_{(n+1)j}^h = 0, \qquad \forall j \in V_c, \forall h \in H, \label{R2} \\
& \sum_{j \in V_c} x_{0j}^0 = 1, \label{R3} \\
& \sum_{j \in V_c} \sum_{h \in H \setminus \{0\}} x_{0j}^h = 0, \label{R4} \\
& \sum_{j \in V \setminus \{0,i\}} \sum_{h \in H} x_{ij}^h = 1, \qquad \forall i \in V_c, \label{R5} \\
& \sum_{i \in V \setminus \{j, n+1\}} \sum_{h \in H} x_{ij}^h = 1, \qquad \forall j \in V_c, \label{R6} \\
& \sum_{i \in V_c} \sum_{h \in H} x_{i(n+1)}^h = 1, \label{R7} \\
& \sum_{j \in V \setminus \{i\}} (y_{ji} - y_{ij})=2, \qquad \forall i \in V_c,\label{R8}\\
& \sum_{j \in V_c} y_{0j}=n, \label{R9} \\
& \sum_{i \in V_c} y_{i0}=0, \label{R10}\\
& \sum_{j \in V_c} y_{(n+1)j}=n, \label{R11}\\
& y_{ij} + y_{ji}= n \sum_{h \in H} (x_{ij}^h+x_{ji}^h), \qquad \forall (i,j) \in A, i < j,\label{R12}\\
& a_0=0, \label{R13} \\
& a_i+s_j+t_{ij}^h-2 \bar{T}|H|(1-x_{ij}^h) \le a_j \label{R14} \\
& a_j \le a_i+s_j+t_{ij}^h+\bar{T}|H|(1-x_{ij}^h), \quad \forall (i,j) \in A, \forall h \in H. \nonumber \\
& \sum_{j \in V \setminus \{0\}} \sum_{h \in H} \bar{T} h x_{ij}^h \le a_i < \sum_{j \in V \setminus \{0\}} \sum_{h \in H} \bar{T} (h+1) x_{ij}^h, \qquad \forall i \in V_c, \label{R15}\\
& a_{n+1} < \bar{T}|H|, \label{R16} \\
& x_{ij}^h \in \{0,1\}, \qquad \forall (i,j) \in A, \forall h \in H, \label{R17} \\
& y_{ij} \ge 0, \qquad \forall (i,j) \in A,\label{R18} \\
& a_{i} \ge 0, \qquad \forall i \in V.\label{R19}  
\end{align}

The objective function \eqref{R0} minimizes the total travel time of the route. Constraints \eqref{R1}--\eqref{R2} ensure that no arc enters the depot and no arc leaves the terminal, so that any feasible solution is a directed path from the depot to the terminal. The departure from the depot is forced to occur in the first time interval  $h \in H$ by constraints \eqref{R3}--\eqref{R4}. The assignment constraints \eqref{R5}--\eqref{R7} then guarantee that the route starts at the depot, serves each customer exactly once, and ends at the terminal. Balance constraints \eqref{R8} at the customer nodes ensure the conservation of the flow associated with customers that have already been served and those that still need to be served. Constraints \eqref{R9}--\eqref{R11} define the corresponding boundary conditions at the depot and at the terminal, while constraints \eqref{R12} link the routing and flow variables, enforce connectivity of the solution, and eliminate subtours. The departure time from the depot at the beginning of the first time interval is fixed by constraint \eqref{R13}. Departure times at the nodes are controlled by constraints \eqref{R14} along the route within the time-discretized framework, while constraints \eqref{R15}--\eqref{R16} identify, for each node, the time interval in which traversing starts and ensure that the route is completed within the planning horizon. Finally, constraints \eqref{R17}--\eqref{R19} specify the domains of the decision variables.


In addition, we include the valid inequality \eqref{R20}, where $L$ denotes a lower bound on the total travel time of any feasible route. The constant $L$ is defined as the sum of the smallest travel times leaving each eligible node over all time intervals, namely, $L= \min_{j \in V_c} t_{0j}^0 + \sum_{i \in V_c} \min_{\substack{j \in V \setminus \{0, i\} \\ h \in H}} t_{ij}^h$. 
This valid inequality is used to further tighten the formulation.

\begin{align}
	\sum_{(i,j) \in A} \sum_{h \in H} x_{ij}^h \ge L.\label{R20}
\end{align}

\subsubsection{Encoding/Decoding}

We now illustrate how our approach applies to the TD-TSP by designing a decoder for the MILP in \eqref{R0}--\eqref{R19}. The decoder enforces integrality and routing constraints by construction: it always generates a simple route that departs from the depot, serves each customer exactly once, and ends at the terminal node $n+1$ (a copy of the depot), so that degree and connectivity requirements are automatically satisfied. The time-dependent structure and the planning horizon are handled by simulating the route over time, with travel times evaluated as a function of the current departure time and the corresponding time interval updated at each step. Whenever the accumulated time exceeds the planning horizon, a large penalty is added to the objective value, discouraging such solutions and guiding the search toward routes that respect the discretized planning horizon.

Each solution is encoded as a random-key vector of length $n$, where $n$ is the number of customers and each key is a real value in $[0,1)$. The decoder maps this vector into a permutation $v=(v_1,\ldots,v_n)$ of $V_c$ by sorting the customers in nondecreasing order of their random-key values. 
This permutation defines how the variables $x^h_{ij}$ are set, with the depot fixed as the starting node $0$ and its copy (the terminal node $n+1$) designated as the end of the route.
An illustrative example is given next.

Given the permutation $v=(v_1,\ldots,v_n)$, the decoder evaluates the time-dependent cost by simulating the route over the discretized planning horizon. The vehicle departs from the depot at time zero and traverses the arc sequence $(0,v_1),(v_1,v_2),\ldots,(v_n,n+1)$, updating the accumulated time and the corresponding departure interval $h \in H$ after each move. Initially, the current node, $c$, is the starting node $0$ and $i=v_1$, then the variable $x^{h}_{c,v_{i}}$ is set as 1. The current customer $c$ is updated to $v_i$, the next customer in the sequence $v$ is set to $v_i$, and the process is repeated. While the accumulated time remains within the planning horizon, the decoder accumulates the corresponding time-dependent travel and service times; once the horizon is exceeded, the visiting order is fixed, and the solution is penalized.

From the MILP perspective, the simulation induces a sparse variable assignment. For each arc $(i,j)$ on the decoded route, the decoder sets $x_{ij}^{h}=1$ for $h=\left\lfloor a_i/\bar{T}\right\rfloor$ and sets all other $x$-variables to zero. The flow variables follow the same incidence pattern: if arc $(i,j)$ is not used, $(y_{ij},y_{ji})=(0,0)$; otherwise, $(y_{ij},y_{ji})$ is assigned to satisfy the flow-balance constraints and \eqref{R12}, implying $y_{ij}+y_{ji}=n$ for each selected pair $\{i,j\}$. Together with the simulated departure times $a_i$, this defines a complete MILP-consistent assignment of the routing, flow, and timing variables for the decoded route.


As in the portfolio optimization example, the random-key representation focuses on the core combinatorial structure of the solution. Here, the random-key vector encodes only the visiting order of the customers, rather than the full set of time-indexed decision variables in the MILP formulation. This substantially reduces the dimensionality of the search space, which is beneficial for the efficiency of the search procedure. In turn, the decoder is responsible for enforcing the routing constraints and for accurately evaluating the time-dependent travel cost and feasibility of each candidate solution. Algorithm \ref{alg:decoder-tdtsp} summarizes the decoder pseudocode, and Figure \ref{fig:exampleDecoder2} presents an example of mapping a random-key vector $\mathcal{X} = \{0.81, 0.32, 0.54, 0.25,$ $0.15, 0.91\}$ for a dataset with $n = 6$ customers. 

\begin{algorithm}[htbp] 
\footnotesize
\DontPrintSemicolon
\KwIn{$n, |H|, \bar{T}, s[i], \forall i \in V, t[i][j][h], \forall (i,j) \in A, \forall h \in H,  \mathcal{X} $}
\KwOut{$x, y, a, cost$}
\ \ $v \leftarrow [1,...n]$\\
Sort $v[1,...,n]$ in nondecreasing order of $\mathcal{X}[v[i]]$\\
$penalty \leftarrow 0$;
$cost \leftarrow 0$; 
$h \leftarrow 0$;
$current \leftarrow 0$;
$flow \leftarrow n$;\\
$a[0] \leftarrow 0$ \\
\For{$i \leftarrow 1$ \KwTo $n$}{
\If{$h < |H|$}{
$x[current][v[i]][h] \leftarrow 1$ \\
$a[v[i]] \leftarrow a[current] + t[current][v[i]][h] + s[v[i]]$ \\
$h \leftarrow \left\lfloor a[v[i]]/\bar{T}\right\rfloor$
}
\Else{
$x[current][v[i]][|H|-1] \leftarrow 1$\\
$a[v[i]] \leftarrow a[current] + t[current][v[i]][|H|-1] + s[v[i]]$
}
$y[current][v[i]]\leftarrow flow$ \\
$y[v[i]][current]\leftarrow n-flow$ \\
$flow \leftarrow flow -1$\\
$current \leftarrow v[i]$
}
\If{$(h < |H|) \wedge ((a[current]+t[current][0][h])<|H|\bar{T})$} {
$x[current][n+1][h] \leftarrow 1$\\
$a[n+1]=a[current]+t[current][0][h]$
    }
    \Else{
    $x[current][n+1][|H|-1] \leftarrow 1$\\
    $a[n+1] \leftarrow a[current] + t[current][0][|H|-1]$\\
    $penalty \leftarrow 10^3$ 
}
\ \ $y[current][n+1]\leftarrow flow$\\
$y[n+1][current]\leftarrow n-flow$\\
$cost \leftarrow \sum_{(i,j) \in A} \sum_{h \in H} t_{ij}^h x_{ij}^h  + |H|\bar{T} \times penalty$ \\
\caption{Decoder for the TD-TSP}
\label{alg:decoder-tdtsp}
\end{algorithm}


\begin{figure}[htbp]
\centering
\begin{adjustbox}{width=\textwidth}
\begin{tikzpicture}[
  font=\small,
  line cap=round,
  line join=round,
  outer/.style={draw, very thick},
  sepLine/.style={draw, thick},
  tbl/.style={
    matrix of nodes,
    nodes={
      draw,
      minimum height=6mm,
      minimum width=10mm,
      align=center,
      inner sep=1.5pt
    },
    row sep=-\pgflinewidth,
    column sep=-\pgflinewidth
  },
  noboxrow/.style={
    matrix of nodes,
    nodes={
      draw=none,
      minimum height=6mm,
      minimum width=10mm,
      align=center,
      inner sep=1.5pt
    },
    row sep=0mm,
    column sep=0mm
  },
  bigtbl/.style={
    matrix of nodes,
    inner sep=0pt,
    outer sep=0pt,
    nodes={
      draw,
      very thin,
      minimum height=5.2mm,
      minimum width=7.0mm,
      align=center,
      inner sep=1.0pt,
      text height=1.8ex,
      text depth=0.6ex
    },
    row sep=-\pgflinewidth,
    column sep=-\pgflinewidth
  }
]

\def\W{20}
\def\H{12.3}

\draw[outer] (0,0) rectangle (\W,\H);

\def\ySplit{6.3}     
\def\yCostSplit{1.5} 

\draw[sepLine] (0,\ySplit) -- (\W,\ySplit);
\draw[sepLine] (0,\yCostSplit) -- (\W,\yCostSplit);

\def\xA{6.6}
\def\xB{13.2}
\def\pad{0.35}

\coordinate (C1NW) at ($(0,\H)+(\pad,-\pad)$);
\coordinate (C2NW) at ($(\xA,\H)+(\pad,-\pad)$);
\coordinate (C3NW) at ($(\xB,\H)+(\pad,-\pad)$);

\node[anchor=north west] (title1) at (C1NW) {\normalsize$n=6;\,|H|=2;\,\bar{T}=30$};

\matrix[noboxrow, anchor=north west] (idx1) at ($(title1.south west)+(0,0)$) {\normalsize
  1 & 2 & 3 & 4 & 5 & 6 \\
};

\matrix[tbl, anchor=north west] (val1) at ($(idx1.south west)+(0,+0.2)$) {\normalsize
  0.81 & 0.32 & 0.54 & 0.29 & 0.15 & 0.91 \\
};

\matrix[noboxrow, anchor=north west] (idx2) at ($(val1.south west)+(0,0)$) {\normalsize
  5 & 4 & 2 & 3 & 1 & 6 \\
};

\matrix[tbl, anchor=north west] (val2) at ($(idx2.south west)+(0,+0.2)$) {\normalsize
  0.15 & 0.29 & 0.32 & 0.54 & 0.81 & 0.91 \\
};

\matrix[noboxrow, anchor=north west] (sidx) at ($(val2.south west)+(0,0)$) {\normalsize
  $s_1$ & $s_2$ & $s_3$ & $s_4$ & $s_5$ & $s_6$ \\
};

\matrix[tbl, anchor=north west] (sval) at ($(sidx.south west)+(0,+0.2)$) {\normalsize
  5 & 5 & 6 & 4 & 3 & 4 \\
};

\node[anchor=north] (h0) at ($(C2NW)+(3.1,0)$) {};
\node[anchor=north west] at (C2NW |- h0.north) {$h=0$};

\matrix[bigtbl, anchor=north west] (M0) at ($(h0.south west)+(-3.1,-0.25)$) {
  $\phantom{0}$ & \textbf{0} & \textbf{1} & \textbf{2} & \textbf{3} & \textbf{4} & \textbf{5} & \textbf{6} & \textbf{7} \\
  \textbf{0}  & 0 & 5 & 7 & 4 & 1 & 3 & 6 & 0 \\
  \textbf{1}  & 4 & 0 & 8 & 1 & 1 & 4 & 2 & 4 \\
  \textbf{2}  & 7 & 8 & 0 & 5 & 2 & 6 & 6 & 7 \\
  \textbf{3}  & 5 & 2 & 4 & 0 & 1 & 3 & 2 & 5 \\
  \textbf{4}  & 3 & 1 & 2 & 1 & 0 & 7 & 8 & 3 \\
  \textbf{5}  & 2 & 3 & 5 & 3 & 9 & 0 & 4 & 2 \\
  \textbf{6}  & 5 & 2 & 8 & 2 & 7 & 2 & 0 & 5 \\
  \textbf{7}  & 0 & 5 & 7 & 4 & 1 & 3 & 6 & 0 \\
};

\node[anchor=north] (h1) at ($(C3NW)+(3.1,0)$) {};
\node[anchor=north west] at (C3NW |- h1.north) {$h=1$};

\matrix[bigtbl, anchor=north west] (M1) at ($(h1.south west)+(-3.1,-0.25)$) {
  $\phantom{0}$ & \textbf{0} & \textbf{1} & \textbf{2} & \textbf{3} & \textbf{4} & \textbf{5} & \textbf{6} & \textbf{7} \\
  \textbf{0}  & 0 & 8 & 10 & 3 & 1 & 2 & 4 & 0 \\
  \textbf{1}  & 7 & 0 & 8  & 1 & 3 & 4 & 4 & 7 \\
  \textbf{2}  & 9 & 8 & 0  & 6 & 2 & 6 & 8 & 9 \\
  \textbf{3}  & 5 & 4 & 4  & 0 & 1 & 3 & 6 & 5 \\
  \textbf{4}  & 2 & 1 & 2  & 1 & 0 & 7 & 8 & 2 \\
  \textbf{5}  & 3 & 2 & 5  & 4 & 11& 0 & 4 & 3 \\
  \textbf{6}  & 5 & 3 & 8  & 7 & 7 & 2 & 0 & 5 \\
  \textbf{7}  & 0 & 8 & 10 & 3 & 1 & 2 & 4 & 0 \\
};

\coordinate (TextNW) at ($(0,\ySplit)+(\pad,-\pad)$);

\pgfmathsetlengthmacro{\TextWidth}{\W cm - 2*\pad cm}

\pgfmathsetlengthmacro{\TextBandHeight}{(\ySplit-\yCostSplit)*1cm - 2*\pad cm}

\node[anchor=north west, inner sep=0pt] at (TextNW) {%
\begin{minipage}[c][\TextBandHeight][c]{\TextWidth}
\large
\raggedright
\(
a[0]=0;\; h=0;\; penalty=0;
\)\par
\(
x[0][5][0]=1;\; y[0][5]=6;\; y[5][0]=0;\;
a[5]=a[0]+t[0][5][0]+s[5]=0+3+3=6;\; h=0;\; penalty=0;
\)\par
\(
x[5][4][0]=1;\; y[5][4]=5;\; y[4][5]=1;\;
a[4]=a[5]+t[5][4][0]+s[4]=6+9+4=19;\; h=0;\; penalty=0;
\)\par
\(
x[4][2][0]=1;\; y[4][2]=4;\; y[2][4]=2;\;
a[2]=a[4]+t[4][2][0]+s[2]=19+2+5=26;\; h=0;\; penalty=0;
\)\par
\(
x[2][3][0]=1;\; y[2][3]=3;\; y[3][2]=3;\;
a[3]=a[2]+t[2][3][0]+s[3]=26+5+6=37;\; h=1;\; penalty=0;
\)\par
\(
x[3][1][1]=1;\; y[3][1]=2;\; y[1][3]=4;\;
a[1]=a[3]+t[3][1][1]+s[1]=37+4+5=46;\; h=1;\; penalty=0;
\)\par
\(
x[1][6][1]=1;\; y[1][6]=1;\; y[6][1]=5;\;
a[6]=a[1]+t[1][6][1]+s[6]=46+4+4=54;\; h=1;\; penalty=0;
\)\par
\(
x[6][7][1]=1;\; y[6][7]=0;\; y[7][6]=6;\;
a[7]=a[6]+t[6][7][1]=54+5=59;\; h=1;\; penalty=0;
\)
\end{minipage}
};

\node[anchor=center, align=center] at (\W/2,\yCostSplit/2) {\large{$
\mathit{Cost}
=
\sum_{(i,j)\in A}\ \sum_{h\in H} t_{ij}^{h}\,x_{ij}^{h}
+|H|\,\bar{T}\,\times penalty
$}};

\end{tikzpicture}
\end{adjustbox}

\caption{Decoder example for the TD-TSP model. The solution is encoded in the random-key vector of size $n$. For the given example, the selected customer visiting sequence is $5-4-2-3-1-6$.}
\label{fig:exampleDecoder2}
\end{figure}

\section{Experiments and Results}
\label{sec:sec5}

This section demonstrates the feasibility and performance of the random-key approach using numerical tests on two sets of benchmark instances. The first set consists of publicly available benchmark instances \footnote{OR-Library: \url{https://www.brunel.ac.uk/~mastjjb/jeb/info.html}} based on real-world data for the limited asset, or cardinality-constrained, Markowitz. The instances comprise expected return vectors and covariance matrices of sizes between 31 and 2151 constructed from weekly priced data for Hang Seng (31 assets), DAX (85 assets), FTSE (89 assets), S\&P (98 assets), Nikkei (225 assets), S\&P 500 (457 assets), Russell 2000 (1318 assets) and Russell 3000 (2151 assets). Each instance was tested with different values of $K$ to analyze the influence of this parameter on the optimization methods. The parameters $l_i$ and $u_i$ were fixed at 0.01 and 0.25, respectively. For the first set of experiments, the mathematical model (\ref{fo})-(\ref{r4}) was solved using Gurobi version 12.0.2 \citep{gurobi} with a time limit of 1800s.

For the TD-TSP experiments, we consider a second set of 225 instances derived from real traffic data. This set covers problems with $n \in \{10,12,14,20,22,$ $24,50,52,54,80,82,84,100,102,104\}$ customers. Each customer $i \in V_c$ is associated with a deterministic service time $s_i$ whose value lies in an interval $[m,M]$ that depends on the size of the instance: $s_i \in [1800,2700]$s for $n \in \{10,12,14\}$, $s_i \in [900,1500]$ for $n \in \{20,22,24\}$, $s_i \in [360,600]$ for $n \in \{50,52,54\}$, $s_i \in [180,360]$ for $n \in \{80,82,84\}$ and $s_i \in [120,240]$ for $n \in \{100,102,104\}$.

The planning horizon is fixed at 54,000s and is discretized into $|H| \in \{5,10,15\}$ time intervals of equal length. For each instance size $n$ we therefore consider three time-discretization schemes, one for each value of $|H|$. For every combination of instance size $n$ and number of time intervals $|H|$, five independent instances are defined, leading to a total of 225 test instances. For this experiment, the mathematical model \eqref{R0}--\eqref{R19}, with the valid inequality \eqref{R20}, was solved using Gurobi version 12.0.2 \citep{gurobi} with a time limit of 1800s. In all TD-TSP experiments, Gurobi was provided with an initial incumbent solution obtained from a greedy heuristic \citep{malandraki1992time}.

The RKO framework was coded in C++ and compiled using GCC. The RKO framework is available at \url{https://github.com/RKO-solver}. All experiments were conducted on a PC with a Dual Xenon Silver 4114 20c/40t 2.2GHz processor, 96GB of DDR4 RAM, and Ubuntu 24.04.2 x64. We employ nine metaheuristics (GA, BRKGA, BRKGA-CS, PSO, SA, LNS, ILS, VNS, and GRASP) and IPR, all of which run in parallel within RKO's internal execution mechanism.

\subsection{Computational Results for the Limited Asset Markowitz}

The RKO was run 30 times for each instance, with different seeds and a CPU time restriction ($\mathit{run_{time}}$). The CPU time limit for the RKO was set in seconds. A step function \(\mathit{run_{time}}(n)\) was constructed to transfer the number of assets to discrete $\mathit{run_{time}}$ levels. The runtime is defined as $run_{time}(n) = v$ for $n$ in the respective intervals, where the pairs $(v, \text{interval})$ are: $(10, n \le 31)$, $(20, 31 < n \le 98)$, $(30, 98 < n \le 225)$, $(50, 225 < n \le 457)$, $(100, 457 < n \le 1318)$, and $(200, n > 1318)$.



Tables \ref{tab:lambda03} -- \ref{tab:lambda07} show the computational results of Gurobi and RKO to solve the portfolio selection problem with three different values of $\lambda=\{0.3, 0.5, 0.7\}$. These tables contain information about the instances ($n$, $K$, $l_i$, $u_i$, and $\lambda$), the results of Gurobi (upper bound ($\mathrm{UB}$), lower bound ($\mathrm{LB}$), $\mathrm{gap(\%)}$, and the total time to find the optimal solution - limited to 1800s), and the results of RKO (objective function of the best solution found in 30 runs ($\mathrm{OFV}$), relative percentage deviation ($\mathrm{RPD}$), and the average computational time to find the best solution). We calculated the $\mathrm{RPD} = (\mathrm{OFV}/\mathrm{UB} - 1) \times 100$, which considers the objective function value of the best solutions found by RKO and the upper bound of Gurobi. We calculated two metrics: first, in terms of the best solution found in 30 runs ($\mathrm{RPD_b}$) and second, as an average of the $\mathrm{RPD}$ found in each of the 30 runs ($\mathrm{RPD_a}$).

\begin{table}[htbp]
\caption{Computational results for Portfolio instances ($\lambda = 0.3$).}
\label{tab:lambda03}
\resizebox{\textwidth}{!}{%
\begin{tabular}{llrrrrrrrrrrrrrr} \toprule
 &  &  &  &  &  &  & \multicolumn{4}{c}{Gurobi} &  & \multicolumn{4}{c}{RKO} \\ \cmidrule{8-11} \cmidrule{13-16}
Instance & $n$ & $K$ & $l_i$ & $u_i$ & $\lambda$ &  & UB & LB & gap(\%) & total time (s) &  & OFV & $\mathrm{RPD_b}$ & $\mathrm{RPD_a}$ & best found at (s) \\ \midrule
port1 (Hang Seng) & 31 & 5 & 0.01 & 0.25 & 0.3 &  & -0.00466 & -0.00466 & 0.00 & 0.02 &  & -0.00466 & 0.00 & 0.01 & 6.64 \\
 &  & 8 &  &  &  &  & -0.00465 & -0.00465 & 0.00 & 0.01 &  & -0.00465 & 0.00 & 0.00 & 9.77 \\
 &  & 10 &  &  &  &  & -0.00464 & -0.00464 & 0.00 & 0.01 &  & -0.00464 & 0.01 & 0.01 & 9.94 \\
 &  & 15 &  &  &  &  & -0.00461 & -0.00461 & 0.00 & 0.01 &  & -0.00461 & 0.03 & 0.10 & 9.97 \\
port2 (DAX) & 85 & 10 & 0.01 & 0.25 & 0.3 &  & -0.00532 & -0.00532 & 0.00 & 0.03 &  & -0.00532 & 0.00 & 0.05 & 19.89 \\
 &  & 20 &  &  &  &  & -0.00515 & -0.00515 & 0.00 & 0.03 &  & -0.00515 & 0.12 & 0.17 & 19.95 \\
 &  & 25 &  &  &  &  & -0.00505 & -0.00505 & 0.00 & 0.03 &  & -0.00504 & 0.29 & 0.38 & 20.00 \\
 &  & 45 &  &  &  &  & -0.00435 & -0.00435 & 0.00 & 0.03 &  & -0.00421 & 3.20 & 3.30 & 19.98 \\
port3 (FTSE) & 89 & 10 & 0.01 & 0.25 & 0.3 &  & -0.00448 & -0.00448 & 0.00 & 0.03 &  & -0.00448 & 0.02 & 0.02 & 19.47 \\
 &  & 20 &  &  &  &  & -0.00439 & -0.00439 & 0.00 & 0.03 &  & -0.00439 & 0.04 & 0.10 & 19.85 \\
 &  & 25 &  &  &  &  & -0.00432 & -0.00432 & 0.00 & 0.03 &  & -0.00432 & 0.05 & 0.10 & 19.95 \\
 &  & 45 &  &  &  &  & -0.00388 & -0.00388 & 0.00 & 0.03 &  & -0.00359 & 7.71 & 7.96 & 20.00 \\
port4 (S\&P) & 98 & 10 & 0.01 & 0.25 & 0.3 &  & -0.00549 & -0.00549 & 0.00 & 0.04 &  & -0.00549 & 0.00 & 0.01 & 19.78 \\
 &  & 20 &  &  &  &  & -0.00531 & -0.00531 & 0.00 & 0.04 &  & -0.00531 & 0.04 & 0.15 & 19.95 \\
 &  & 30 &  &  &  &  & -0.00507 & -0.00507 & 0.00 & 0.04 &  & -0.00505 & 0.40 & 0.53 & 19.99 \\
 &  & 50 &  &  &  &  & -0.00440 & -0.00440 & 0.00 & 0.04 &  & -0.00407 & 7.48 & 7.52 & 19.99 \\
port5 (Nikkei) & 225 & 10 & 0.01 & 0.25 & 0.3 &  & -0.00230 & -0.00230 & 0.00 & 0.12 &  & -0.00230 & 0.02 & 0.02 & 28.45 \\
 &  & 20 &  &  &  &  & -0.00219 & -0.00219 & 0.00 & 0.12 &  & -0.00219 & 0.03 & 0.51 & 29.95 \\
 &  & 30 &  &  &  &  & -0.00202 & -0.00202 & 0.00 & 0.12 &  & -0.00198 & 2.40 & 3.54 & 29.96 \\
 &  & 50 &  &  &  &  & -0.00158 & -0.00158 & 0.00 & 0.12 &  & -0.00143 & 9.51 & 10.81 & 29.92 \\
port6 (S\&P 500) & 457 & 20 & 0.01 & 0.25 & 0.3 &  & -0.00888 & -0.00888 & 0.00 & 338.80 &  & -0.00888 & 0.03 & 0.04 & 49.94 \\
 &  & 40 &  &  &  &  & -0.00836 & -0.00839 & 0.40 & 1800.99 &  & -0.00815 & 2.46 & 2.74 & 49.84 \\
 &  & 50 &  &  &  &  & -0.00799 & -0.00807 & 1.00 & 1800.63 &  & -0.00749 & 6.29 & 7.48 & 49.92 \\
 &  & 70 &  &  &  &  & -0.00711 & -0.00729 & 2.48 & 1800.99 &  & -0.00651 & 8.46 & 10.58 & 49.94 \\
port7 (Russell 2000) & 1318 & 30 & 0.01 & 0.25 & 0.3 &  & -0.01226 & -0.01239 & 1.05 & 1808.74 &  & -0.01224 & 0.16 & 0.38 & 99.87 \\
 &  & 50 &  &  &  &  & -0.01165 & -0.01190 & 2.14 & 1812.12 &  & -0.01151 & 1.24 & 1.41 & 99.98 \\
 &  & 60 &  &  &  &  & -0.01126 & -0.01157 & 2.76 & 1813.51 &  & -0.01115 & 0.99 & 1.07 & 99.96 \\
 &  & 80 &  &  &  &  & -0.01034 & -0.01080 & 4.48 & 1805.09 &  & -0.01019 & 1.44 & 2.95 & 99.51 \\
port8 (Russell 3000) & 2151 & 40 & 0.01 & 0.25 & 0.3 &  & -0.00725 & -0.00732 & 0.86 & 1814.92 &  & -0.00709 & 2.27 & 2.31 & 199.67 \\
 &  & 60 &  &  &  &  & -0.00624 & -0.00689 & 10.49 & 1810.55 &  & -0.00643 & -3.11 & -2.79 & 199.86 \\
 &  & 70 &  &  &  &  & -0.00588 & -0.00669 & 13.78 & 1810.11 &  & -0.00616 & -4.63 & -4.26 & 199.96 \\
 &  & 90 &  &  &  &  & -0.00534 & -0.00609 & 14.00 & 1810.22 &  & -0.00542 & -1.49 & -1.04 & 199.86
 \\ \bottomrule
\end{tabular}%
}
\end{table}

\begin{table}[htbp]
\caption{Computational results for Portfolio instances ($\lambda = 0.5$).}
\label{tab:lambda05}
\resizebox{\textwidth}{!}{%
\begin{tabular}{llrrrrrrrrrrrrrr} \toprule
 &  &  &  &  &  &  & \multicolumn{4}{c}{Gurobi} &  & \multicolumn{4}{c}{RKO} \\ \cmidrule{8-11} \cmidrule{13-16}
Instance & $n$ & $K$ & $l_i$ & $u_i$ & $\lambda$ &  & UB & LB & gap(\%) & total time (s) &  & OFV & $\mathrm{RPD_b}$ & $\mathrm{RPD_a}$ & best found at (s) \\ \midrule
port1 (Hang Seng) & 31 & 5 & 0.01 & 0.25 & 0.5 &  & -0.00297 & -0.00297 & 0.00 & 0.43 &  & -0.00297 & 0.00 & 0.00 & 8.36 \\
 &  & 8 &  &  &  &  & -0.00297 & -0.00297 & 0.00 & 0.03 &  & -0.00297 & 0.02 & 0.02 & 9.93 \\
 &  & 10 &  &  &  &  & -0.00296 & -0.00296 & 0.00 & 0.01 &  & -0.00296 & 0.00 & 0.01 & 9.91 \\
 &  & 15 &  &  &  &  & -0.00293 & -0.00293 & 0.00 & 0.01 &  & -0.00293 & 0.08 & 0.15 & 9.98 \\
port2 (DAX) & 85 & 10 & 0.01 & 0.25 & 0.5 &  & -0.00365 & -0.00365 & 0.00 & 0.03 &  & -0.00365 & 0.00 & 0.03 & 19.22 \\
 &  & 20 &  &  &  &  & -0.00354 & -0.00354 & 0.00 & 0.03 &  & -0.00354 & 0.00 & 0.16 & 19.79 \\
 &  & 25 &  &  &  &  & -0.00347 & -0.00347 & 0.00 & 0.03 &  & -0.00346 & 0.26 & 0.33 & 20.00 \\
 &  & 45 &  &  &  &  & -0.00298 & -0.00298 & 0.00 & 0.03 &  & -0.00290 & 2.85 & 2.93 & 19.95 \\
port3 (FTSE) & 89 & 10 & 0.01 & 0.25 & 0.5 &  & -0.00299 & -0.00299 & 0.00 & 0.05 &  & -0.00299 & 0.02 & 0.02 & 19.83 \\
 &  & 20 &  &  &  &  & -0.00293 & -0.00293 & 0.00 & 0.03 &  & -0.00293 & 0.00 & 0.03 & 19.84 \\
 &  & 25 &  &  &  &  & -0.00288 & -0.00288 & 0.00 & 0.03 &  & -0.00288 & 0.02 & 0.06 & 19.82 \\
 &  & 45 &  &  &  &  & -0.00260 & -0.00260 & 0.00 & 0.03 &  & -0.00246 & 5.56 & 7.04 & 19.96 \\
port4 (S\&P) & 98 & 10 & 0.01 & 0.25 & 0.5 &  & -0.00359 & -0.00359 & 0.00 & 0.04 &  & -0.00359 & 0.00 & 0.00 & 19.94 \\
 &  & 20 &  &  &  &  & -0.00347 & -0.00347 & 0.00 & 0.04 &  & -0.00347 & 0.22 & 0.30 & 19.92 \\
 &  & 30 &  &  &  &  & -0.00333 & -0.00333 & 0.00 & 0.04 &  & -0.00330 & 0.72 & 0.72 & 19.97 \\
 &  & 50 &  &  &  &  & -0.00294 & -0.00294 & 0.00 & 0.04 &  & -0.00278 & 5.40 & 5.71 & 19.98 \\
port5 (Nikkei) & 225 & 10 & 0.01 & 0.25 & 0.5 &  & -0.00144 & -0.00144 & 0.00 & 0.19 &  & -0.00144 & 0.00 & 0.00 & 29.81 \\
 &  & 20 &  &  &  &  & -0.00136 & -0.00136 & 0.00 & 0.12 &  & -0.00136 & 0.08 & 0.68 & 29.98 \\
 &  & 30 &  &  &  &  & -0.00124 & -0.00124 & 0.00 & 0.13 &  & -0.00123 & 0.99 & 2.73 & 29.85 \\
 &  & 50 &  &  &  &  & -0.00092 & -0.00092 & 0.00 & 0.12 &  & -0.00085 & 8.06 & 10.01 & 29.91 \\
port6 (S\&P 500) & 457 & 20 & 0.01 & 0.25 & 0.5 &  & -0.00510 & -0.00512 & 0.49 & 1800.86 &  & -0.00510 & 0.02 & 0.02 & 49.74 \\
 &  & 40 &  &  &  &  & -0.00488 & -0.00500 & 2.38 & 1800.56 &  & -0.00484 & 0.98 & 1.01 & 49.99 \\
 &  & 50 &  &  &  &  & -0.00472 & -0.00488 & 3.29 & 1800.52 &  & -0.00457 & 3.34 & 3.98 & 49.90 \\
 &  & 70 &  &  &  &  & -0.00431 & -0.00458 & 6.23 & 1800.53 &  & -0.00404 & 6.22 & 8.24 & 49.89 \\
port7 (Russell 2000) & 1318 & 30 & 0.01 & 0.25 & 0.5 &  & -0.00768 & -0.00809 & 5.30 & 1811.59 &  & -0.00776 & -1.04 & -0.98 & 99.69 \\
 &  & 50 &  &  &  &  & -0.00723 & -0.00804 & 11.22 & 1805.69 &  & -0.00740 & -2.35 & -1.78 & 99.97 \\
 &  & 60 &  &  &  &  & -0.00685 & -0.00791 & 15.43 & 1803.82 &  & -0.00714 & -4.24 & -4.18 & 99.87 \\
 &  & 80 &  &  &  &  & -0.00623 & -0.00747 & 19.91 & 1805.34 &  & -0.00657 & -5.56 & -4.92 & 99.44 \\
port8 (Russell 3000) & 2151 & 40 & 0.01 & 0.25 & 0.5 &  & -0.00141 & -2.45691 & 174535.75 & 1810.29 &  & -0.00465 & -230.59 & -230.35 & 199.81 \\
 &  & 60 &  &  &  &  & -0.00111 & -0.00631 & 470.25 & 1809.97 &  & -0.00425 & -283.99 & -282.82 & 199.95 \\
 &  & 70 &  &  &  &  & -0.00087 & -0.00631 & 626.67 & 1810.02 &  & -0.00408 & -370.02 & -367.30 & 199.48 \\
 &  & 90 &  &  &  &  & -0.00074 & -0.00628 & 745.15 & 1810.19 &  & -0.00362 & -387.55 & -385.91 & 199.79
 \\ \bottomrule
\end{tabular}%
}
\end{table}


\begin{table}[htbp]
\caption{Computational results for Portfolio instances ($\lambda = 0.7$)}
\label{tab:lambda07}
\resizebox{\textwidth}{!}{%
\begin{tabular}{llrrrrrrrrrrrrrr} \toprule
 &  &  &  &  &  &  & \multicolumn{4}{c}{Gurobi} &  & \multicolumn{4}{c}{RKO} \\ \cmidrule{8-11} \cmidrule{13-16}
Instance & $n$ & $K$ & $l_i$ & $u_i$ & $\lambda$ &  & UB & LB & gap(\%) & total time (s) &  & OFV & $\mathrm{RPD_b}$ & $\mathrm{RPD_a}$ & best found at (s) \\ \midrule
port1 (Hang Seng) & 31 & 5 & 0.01 & 0.25 & 0.7 &  & -0.00130 & -0.00130 & 0.00 & 0.01 &  & -0.00130 & 0.00 & 0.00 & 9.81 \\
 &  & 8 &  &  &  &  & -0.00130 & -0.00130 & 0.00 & 0.01 &  & -0.00130 & 0.00 & 0.02 & 9.82 \\
 &  & 10 &  &  &  &  & -0.00129 & -0.00129 & 0.00 & 0.01 &  & -0.00129 & 0.00 & 0.03 & 9.23 \\
 &  & 15 &  &  &  &  & -0.00127 & -0.00127 & 0.00 & 0.01 &  & -0.00127 & 0.02 & 0.02 & 9.96 \\
port2 (DAX) & 85 & 10 & 0.01 & 0.25 & 0.7 &  & -0.00198 & -0.00198 & 0.00 & 0.03 &  & -0.00198 & 0.01 & 0.01 & 19.81 \\
 &  & 20 &  &  &  &  & -0.00193 & -0.00193 & 0.00 & 0.03 &  & -0.00192 & 0.06 & 0.12 & 19.78 \\
 &  & 25 &  &  &  &  & -0.00188 & -0.00188 & 0.00 & 0.03 &  & -0.00188 & 0.13 & 0.27 & 19.98 \\
 &  & 45 &  &  &  &  & -0.00162 & -0.00162 & 0.00 & 0.03 &  & -0.00155 & 4.41 & 4.56 & 19.95 \\
port3 (FTSE) & 89 & 10 & 0.01 & 0.25 & 0.7 &  & -0.00155 & -0.00155 & 0.00 & 0.03 &  & -0.00155 & 0.02 & 0.02 & 19.93 \\
 &  & 20 &  &  &  &  & -0.00153 & -0.00153 & 0.00 & 0.03 &  & -0.00153 & 0.00 & 0.03 & 19.87 \\
 &  & 25 &  &  &  &  & -0.00150 & -0.00150 & 0.00 & 0.04 &  & -0.00150 & 0.00 & 0.04 & 19.97 \\
 &  & 45 &  &  &  &  & -0.00136 & -0.00136 & 0.00 & 0.03 &  & -0.00127 & 7.16 & 7.54 & 19.97 \\
port4 (S\&P) & 98 & 10 & 0.01 & 0.25 & 0.7 &  & -0.00175 & -0.00175 & 0.00 & 0.04 &  & -0.00175 & 0.00 & 0.00 & 19.83 \\
 &  & 20 &  &  &  &  & -0.00172 & -0.00172 & 0.00 & 0.04 &  & -0.00172 & 0.07 & 0.12 & 19.83 \\
 &  & 30 &  &  &  &  & -0.00167 & -0.00167 & 0.00 & 0.04 &  & -0.00166 & 0.65 & 0.66 & 19.96 \\
 &  & 50 &  &  &  &  & -0.00152 & -0.00152 & 0.00 & 0.04 &  & -0.00149 & 2.59 & 2.70 & 19.93 \\
port5 (Nikkei) & 225 & 10 & 0.01 & 0.25 & 0.7 &  & -0.00059 & -0.00059 & 0.00 & 0.12 &  & -0.00059 & 0.00 & 0.00 & 28.03 \\
 &  & 20 &  &  &  &  & -0.00056 & -0.00056 & 0.00 & 0.13 &  & -0.00056 & 0.25 & 0.32 & 29.96 \\
 &  & 30 &  &  &  &  & -0.00050 & -0.00050 & 0.00 & 0.13 &  & -0.00048 & 3.17 & 4.25 & 29.91 \\
 &  & 50 &  &  &  &  & -0.00032 & -0.00032 & 0.00 & 0.12 &  & -0.00029 & 11.40 & 17.85 & 29.95 \\
port6 (S\&P 500) & 457 & 20 & 0.01 & 0.25 & 0.7 &  & -0.00222 & -0.00243 & 9.60 & 1800.64 &  & -0.00225 & -1.41 & -1.41 & 49.73 \\
 &  & 40 &  &  &  &  & -0.00216 & -0.00244 & 13.11 & 1800.56 &  & -0.00217 & -0.47 & -0.47 & 49.88 \\
 &  & 50 &  &  &  &  & -0.00205 & -0.00239 & 16.89 & 1800.81 &  & -0.00211 & -3.20 & -2.54 & 49.92 \\
 &  & 70 &  &  &  &  & -0.00190 & -0.00239 & 25.41 & 1800.72 &  & -0.00189 & 0.83 & 2.14 & 49.97 \\
port7 (Russell 2000) & 1318 & 30 & 0.01 & 0.25 & 0.7 &  & -0.00258 & -0.00665 & 158.00 & 1803.77 &  & -0.00385 & -49.51 & -49.51 & 99.64 \\
 &  & 50 &  &  &  &  & -0.00232 & -0.00668 & 187.56 & 1803.54 &  & -0.00377 & -62.17 & -61.87 & 99.94 \\
 &  & 60 &  &  &  &  & -0.00101 & -0.00654 & 544.84 & 1803.48 &  & -0.00367 & -261.64 & -261.48 & 99.85 \\
 &  & 80 &  &  &  &  & -0.00092 & -0.00670 & 625.06 & 1803.53 &  & -0.00339 & -266.60 & -265.23 & 99.93 \\
port8 (Russell 3000) & 2151 & 40 & 0.01 & 0.25 & 0.7 &  & -0.00107 & -0.00786 & 637.80 & 1810.10 &  & -0.00235 & -120.88 & -120.82 & 199.84 \\
 &  & 60 &  &  &  &  & -0.00046 & -0.00786 & 1605.93 & 1820.00 &  & -0.00220 & -376.54 & -376.36 & 199.95 \\
 &  & 70 &  &  &  &  & -0.00020 & -0.00787 & 3828.02 & 1810.07 &  & -0.00214 & -966.54 & -960.25 & 199.87 \\
 &  & 90 &  &  &  &  & -0.00035 & -0.00787 & 2161.67 & 1810.22 &  & -0.00192 & -450.98 & -449.60 & 199.96
 \\ \bottomrule
\end{tabular}%
}
\end{table}


Looking at the tables, we note that Gurobi solved up to the Nikkei instance in a few seconds. However, starting with the S\&P 500 instances, Gurobi solved only one instance within the 1800-second limit (S\&P 500 with $K=20$ and $\lambda=0.3$). As the value of $\lambda$ increases, the gaps between instances become larger. This behavior shows that instances become more difficult as we increase the significance of minimizing risk. On the other hand, RKO does not show significant differences for the values of $\lambda$. Considering the first five sets of instances (solved by Gurobi), the average $RPD_b$ for $\lambda=0.3,0.5,$ and $0.7$ was 1.57, 1.21, and 1.50, respectively. However, RKO, like Gurobi, has more difficulty finding good solutions as we increase the value of $K$. We can observe that RKO finds the optimal solution or solutions close to the optimal for the instances solved by Gurobi, except for the instances with the highest values of $K$. For instances in which Gurobi did not find the optimal solution in 1800s, RKO produces solutions close to the upper bound for the S\&P 500 set in 50s. For the Russell 2000 set, the RKO solutions were close to the upper bound in instances with $\lambda=0.3$, and for instances with $\lambda=0.5$ and $0.7$, the results of RKO running in 100s were better than Gurobi's upper bound (except for $\lambda=0.7$ and $K=70$). For the set with the largest number of assets (Russell 3000), RKO finds solutions much better than Gurobi's upper bound in all instances (except for $\lambda=0.3$ and $K=40$). Especially for instances with values of $\lambda=0.5$ and $0.7$, where Gurobi does not find good bounds in 1800s, RKO's results were significantly better in 200s of execution. Overall, Gurobi solved 61 of 96 instances and found the best upper bounds in 13 instances, while RKO found the best upper bound in 40 instances, 22 of which were better than Gurobi's upper bound.

\begin{figure}[htbp]
    \centering
    \includegraphics[width=0.79\linewidth]{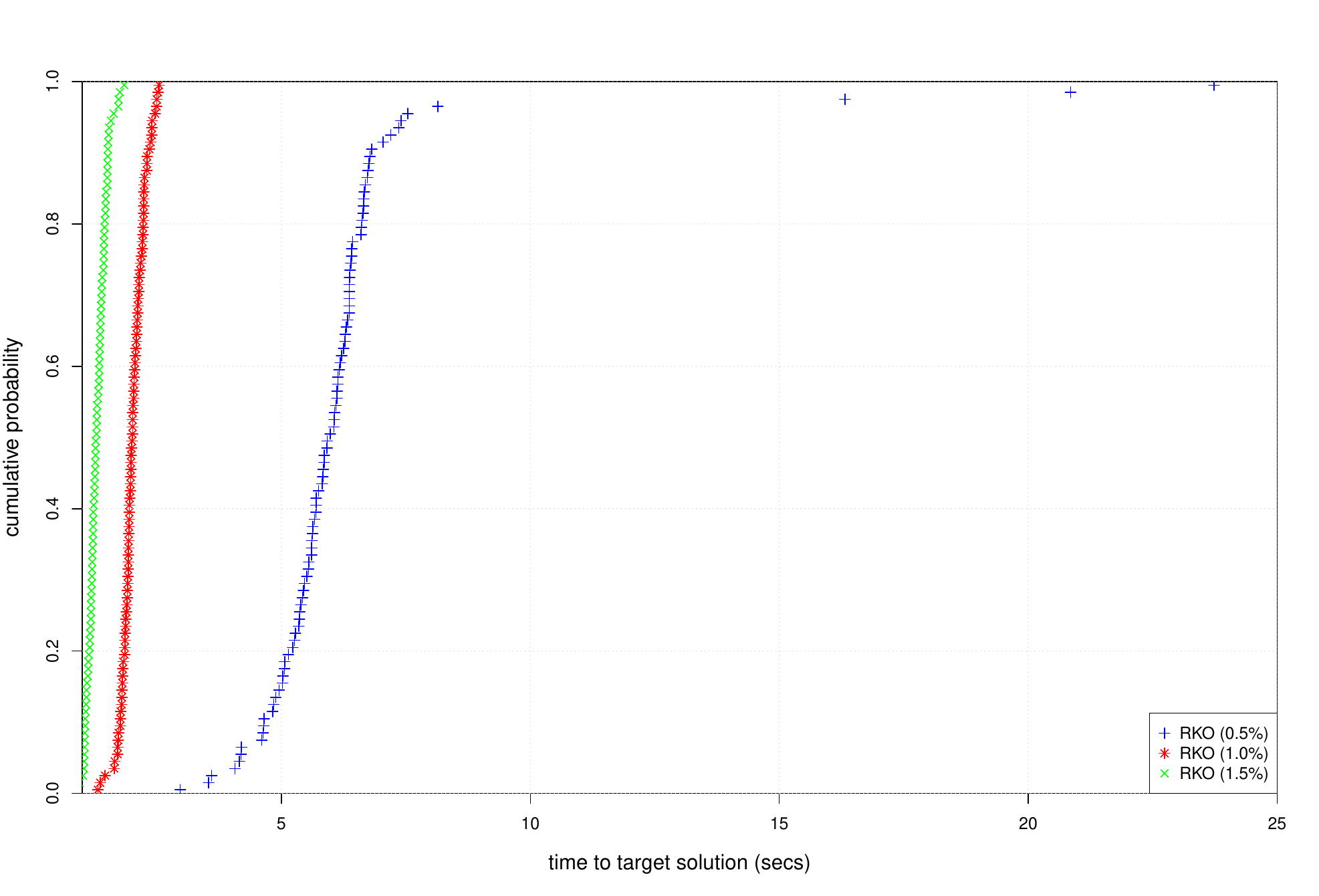}
    \caption{Time-to-target plot for Russell 3000 instance with $K=40$ and $\lambda=0.5$.}
    \label{fig:tttPort8-40}
\end{figure}
\begin{figure}[htbp]
    \centering
    \includegraphics[width=0.79\linewidth]{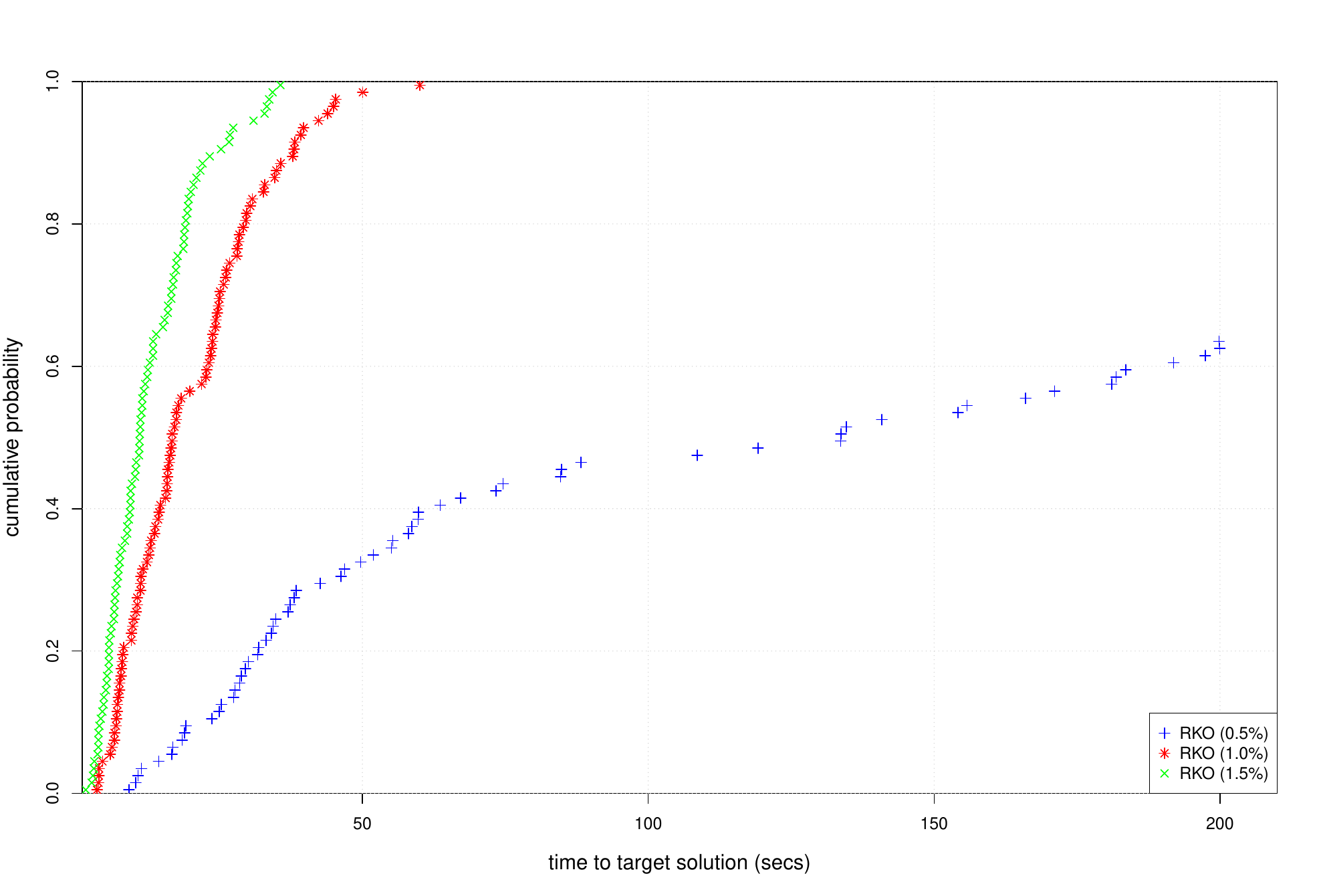}
    \caption{Time-to-target plot for Russell 3000 instance with $K=60$ and $\lambda=0.5$.}
    \label{fig:tttPort8-60}
\end{figure}

Figures \ref{fig:tttPort8-40} and \ref{fig:tttPort8-60} present time-to-target (TTT) plots for the Russell 3000 instances with $\lambda = 0.5$ and $K = 40$ and $60$. The target solution was defined as a percentage deviation from the best-known upper bound, using three thresholds: 0.5\%, 1.0\%, and 1.5\%. The RKO algorithm was executed 100 times, stopping upon finding a solution at least as good as the target or reaching the predefined time limit.
For example, with $K = 40$, RKO had a 96\% probability of finding a solution within 0.5\% of the best known bound in under eight seconds. For the 1.0\% and 1.5\% thresholds, RKO reached the target in 100\% of the runs, within 2.5s and 1.8s, respectively.
In contrast, for $K = 60$, RKO achieved the 0.5\% target in only 65\% of the runs. The probability of reaching this target was 46\% within 100s and 34\% within 50s. However, when the target was relaxed to 1.0\% and 1.5\%, RKO consistently found a solution at the target in under 50s. For instance, with an 80\% probability, the algorithm found the 1.0\% target in 32s and the 1.5\% target in 19s. On average, RKO reached the 1.5\% target in just 64\% of the time required to achieve the 1.0\% target at the same confidence level.

\vspace{-0.5cm}

\vspace{-0.5cm}

\begin{figure}[htbp]
    \centering

    \begin{subfigure}[b]{0.45\textwidth}
        \includegraphics[width=\linewidth]{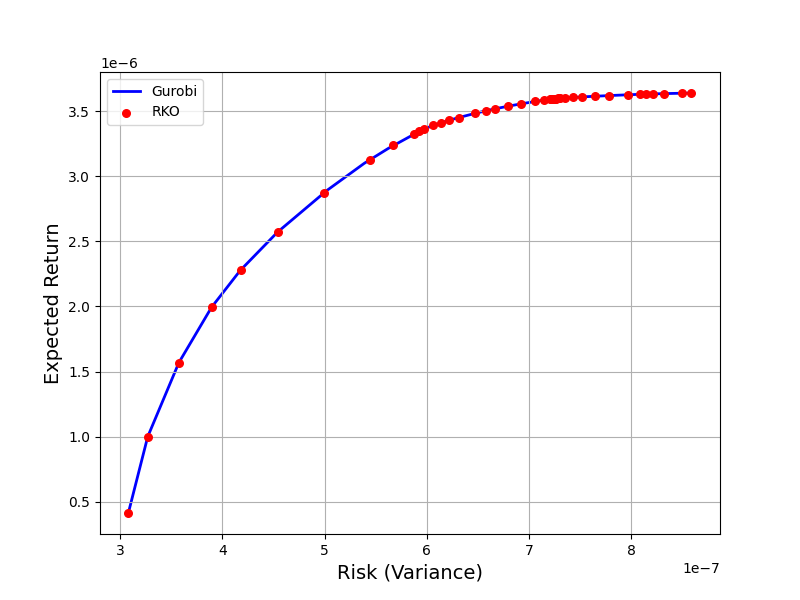}
        \caption{$K=10$}
        \label{fig:Pareto10}
    \end{subfigure}
    \hfill
    \begin{subfigure}[b]{0.45\textwidth}
        \includegraphics[width=\linewidth]{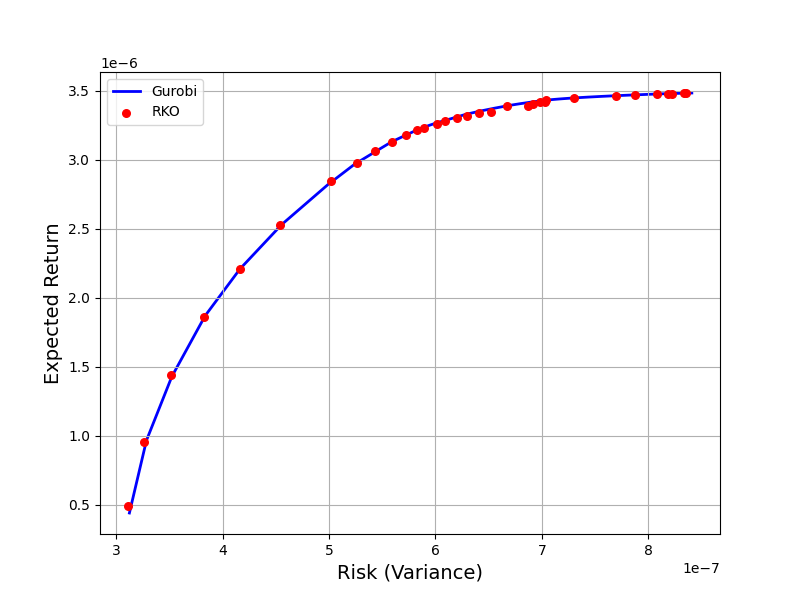}
        \caption{$K=20$}
        \label{fig:Pareto20}
    \end{subfigure}

    \vspace{0.5cm}

    \begin{subfigure}[b]{0.45\textwidth}
        \includegraphics[width=\linewidth]{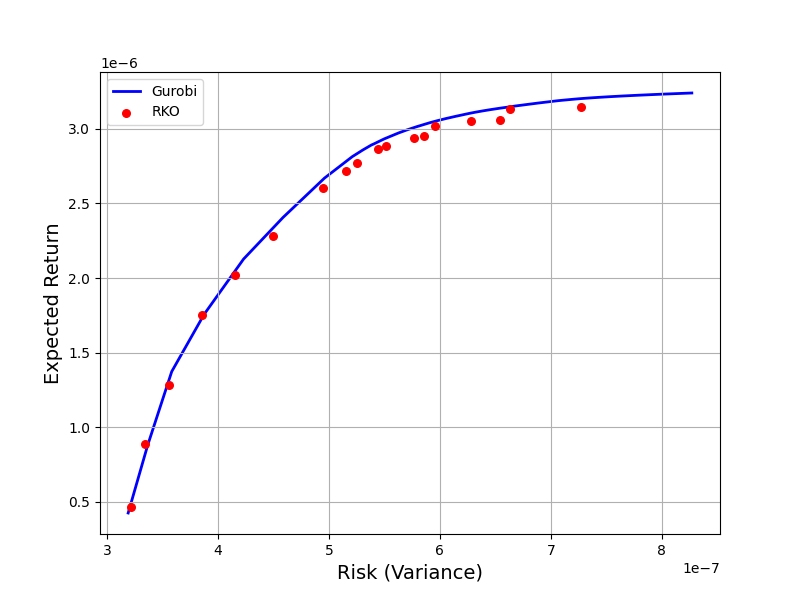}
        \caption{$K=30$}
        \label{fig:Pareto30}
    \end{subfigure}
    \hfill
    \begin{subfigure}[b]{0.45\textwidth}
        \includegraphics[width=\linewidth]{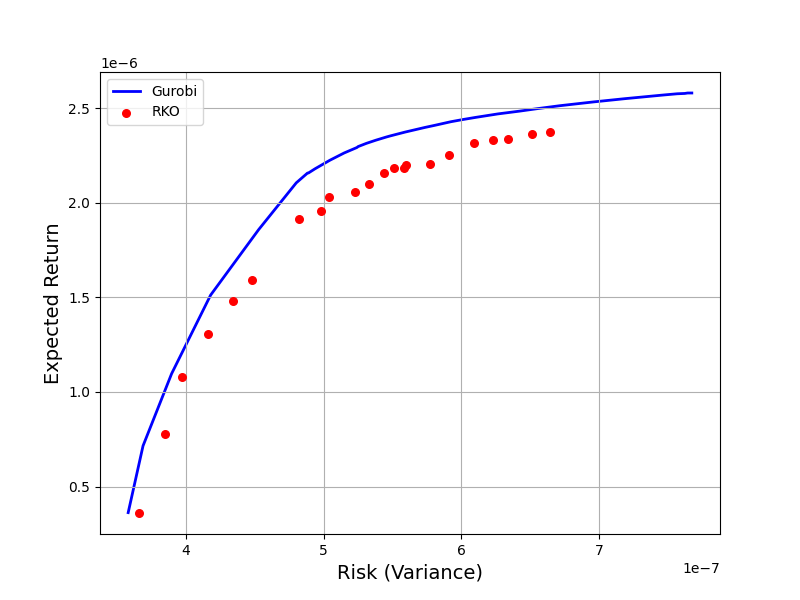}
        \caption{$K=50$}
        \label{fig:Pareto50}
    \end{subfigure}

    \caption{Pareto Frontiers for the Nikkei instances.}
    \label{fig:ParetoNikkei}
\end{figure}

The efficient frontier is the set of portfolios that either maximizes the expected return for a specific level of risk or minimizes the risk for a given expected return. Therefore, the portfolios' position along the efficient frontier indicates the risk-return trade-off.
We generated efficient frontiers for Nikkei instances with $K=10, 20, 30$, and $50$. These instances were selected because they were the largest instances solved by Gurobi. RKO and Gurobi were run for each  Nikkei instance with a sequence of values of $\lambda = \{0.02, 0.04,..., 0.96, 0.98\}$. Figures \ref{fig:Pareto10} - \ref{fig:Pareto50} show the solutions found by Gurobi and RKO. In these figures, the efficient frontiers of Gurobi are represented by a blue line, while the red dots indicate the solutions generated by RKO.
Based on this, we can note that RKO finds solutions on the Pareto frontier for the values of $K=10$ and $K=20$. With $K=30$, RKO finds efficient solutions for values of $\lambda$ that prioritize risk minimization, and as $\lambda$ increases, the solutions move away from the frontier. This behavior was also observed in Tables \ref{tab:lambda03} to \ref{tab:lambda07}. For $K=50$, RKO did not find any efficient solutions on the frontier.

\subsection{Computational Results for the TD-TSP}

Tables \ref{tab:Small} -- \ref{tab:Large} report the average computational results obtained by Gurobi and by RKO for the TD-TSP with three different values of $|H| \in \{5,10,15\}$. For each combination of instance customers $n$ and number of time intervals $|H|$, the results are averaged over the corresponding set of five instances. The CPU time limit for the RKO was set in $n$ seconds. The tables include information on the instances (number of customers $n$), the results of Gurobi (upper bound ($\mathrm{UB}$), lower bound ($\mathrm{LB}$), optimality $\mathrm{gap (\%)}$, and total time required to find an optimal solution, limited to 1800s), and the results of RKO (objective function value of the best solution found in five runs ($\mathrm{OFV}$), relative percentage deviation ($\mathrm{RPD}$), and average computational time to obtain the best solution). Again, we calculated the $\mathrm{RPD} = (\mathrm{OFV}/\mathrm{UB} - 1) \times 100$. Two metrics are reported: the $\mathrm{RPD}$ of the best solution over five runs ($\mathrm{RPD_b}$) and the average $\mathrm{RPD}$ over the five runs ($\mathrm{RPD_a}$). The number of runs was set to five for the TD-TSP to balance statistical significance with the total computational time required.

\begin{table}[htbp]
	\caption{Average computational results for TD-TSP instances with $|H|=5$.}
	\label{tab:Small}
	\resizebox{\textwidth}{!}{%
		\begin{tabular}{lrrrrrrrrrr} \toprule
			&  &\multicolumn{4}{c}{Gurobi} &  & \multicolumn{4}{c}{RKO} \\ \cmidrule{3-6} \cmidrule{8-11}
			$n$ & & UB & LB & gap(\%) & total time (s) &  & OFV & $\mathrm{RPD_b}$ & $\mathrm{RPD_a}$ & best found at (s) \\ \midrule 
			10	&&	7235.60	&	7235.60	&	0.00	&	2.70	&&	7235.60	&	0.00	&	0.00	&	0.00	\\
12	&&	8210.20	&	8210.20	&	0.00	&	14.66	&&	8210.20	&	0.00	&	0.00	&	0.01	\\
14	&&	8930.80	&	8930.80	&	0.00	&	75.29	&&	8930.80	&	0.00	&	0.00	&	0.01	\\
20	&&	11006.00	&	10349.40	&	6.26	&	1618.33	&&	10848.60	&	-0.89	&	-0.89	&	0.47	\\
22	&&	12393.20	&	10887.00	&	12.01	&	1800.05	&&	11877.20	&	-1.98	&	-1.98	&	0.55	\\
24	&&	12479.20	&	10608.60	&	14.67	&	1800.05	&&	11596.80	&	-2.80	&	-2.80	&	1.15	\\
50	&&	20605.40	&	13206.20	&	35.76	&	1801.84	&&	17919.40	&	-11.94	&	-10.52	&	21.55	\\
52	&&	20864.20	&	13820.40	&	33.73	&	1801.55	&&	18097.00	&	-12.16	&	-10.08	&	24.02	\\
54	&&	21064.40	&	13595.00	&	35.39	&	1800.52	&&	18670.00	&	-9.90	&	-7.73	&	26.57	\\
80	&&	26290.40	&	17663.80	&	32.86	&	1800.10	&&	23898.20	&	-8.62	&	-5.66	&	42.40	\\
82	&&	24154.00	&	15185.00	&	37.14	&	1800.07	&&	22337.40	&	-7.17	&	-4.43	&	52.06	\\
84	&&	25423.80	&	16707.00	&	34.33	&	1800.11	&&	22926.00	&	-9.27	&	-5.77	&	55.04	\\
100	&&	24573.80	&	16524.40	&	32.38	&	1800.36	&&	23199.20	&	-5.59	&	-1.52	&	56.38	\\
102	&&	28197.80	&	18692.60	&	33.83	&	1803.57	&&	26560.40	&	-5.78	&	-1.47	&	71.19	\\
104	&&	25295.20	&	16991.60	&	32.30	&	1804.12	&&	24087.80	&	-4.14	&	0.15	&	75.64	\\
\bottomrule
Average &&	18448.27	&	13240.51	&	22.71	&	1434.89	&&	17092.97	&	-5.35	&	-3.51	&	28.47\\
\bottomrule
\end{tabular}%
}
\end{table}

\begin{table}[htbp]
	\caption{Average computational results for TD-TSP instances with $|H|=10$.}
	\label{tab:Medium}
	\resizebox{\textwidth}{!}{%
		\begin{tabular}{lrrrrrrrrrr} \toprule
			&  &\multicolumn{4}{c}{Gurobi} &  & \multicolumn{4}{c}{RKO} \\ \cmidrule{3-6} \cmidrule{8-11}
			$n$ & & UB & LB & gap(\%) & total time (s) &  & OFV & $\mathrm{RPD_b}$ & $\mathrm{RPD_a}$ & best found at (s) \\ \midrule
			10	&&	7354.40	&	7354.40	&	0.00	&	18.32	&&	7354.40	&	0.00	&	0.00	&	0.00	\\
12	&&	8239.20	&	8239.20	&	0.00	&	66.54	&&	8239.20	&	0.00	&	0.00	&	0.01	\\
14	&&	9007.00	&	8828.40	&	1.68	&	590.47	&&	8987.80	&	0.00	&	0.00	&	0.02	\\
20	&&	12339.60	&	9759.60	&	20.64	&	1800.20	&&	11052.80	&	-3.04	&	-3.04	&	0.22	\\
22	&&	12535.80	&	10367.80	&	17.20	&	1800.08	&&	12123.60	&	-3.20	&	-3.20	&	1.21	\\
24	&&	12524.40	&	9966.40	&	20.73	&	1800.10	&&	11877.40	&	-5.00	&	-5.00	&	3.00	\\
50	&&	21231.60	&	12506.20	&	40.99	&	1800.16	&&	18164.60	&	-13.86	&	-12.66	&	21.01	\\
52	&&	21470.20	&	13107.40	&	38.97	&	1800.72	&&	18396.20	&	-13.76	&	-12.30	&	23.30	\\
54	&&	21912.60	&	13113.60	&	39.79	&	1800.08	&&	18599.40	&	-13.34	&	-10.67	&	29.75	\\
80	&&	27278.00	&	16818.60	&	38.49	&	1801.23	&&	24441.40	&	-10.29	&	-8.14	&	48.99	\\
82	&&	24654.00	&	14638.40	&	40.61	&	1805.77	&&	22820.20	&	-7.38	&	-4.37	&	39.78	\\
84	&&	26589.20	&	15607.40	&	41.27	&	1807.88	&&	23150.00	&	-12.62	&	-9.13	&	51.20	\\
100	&&	26349.80	&	15551.40	&	40.11	&	1819.07	&&	23676.40	&	-8.99	&	-4.71	&	72.59	\\
102	&&	29781.60	&	17281.00	&	42.08	&	1800.42	&&	27155.80	&	-8.29	&	-4.93	&	72.60	\\
104	&&	26554.60	&	15695.80	&	41.05	&	1800.45	&&	24095.40	&	-7.91	&	-4.15	&	70.38	\\
			\bottomrule
			Average	&&	19188.13	&	12589.04	&	28.24	&	1487.43	&&	17342.31	&	-7.18	&	-5.49	&	28.94\\
			\bottomrule
		\end{tabular}%
	}
\end{table}

\begin{table}[htbp]
	\caption{Average computational results for TD-TSP instances with $|H|=15$.}
	\label{tab:Large}
	\resizebox{\textwidth}{!}{%
		\begin{tabular}{lrrrrrrrrrr} \toprule
			&  &\multicolumn{4}{c}{Gurobi} &  & \multicolumn{4}{c}{RKO} \\ \cmidrule{3-6} \cmidrule{8-11}
			$n$ & & UB & LB & gap(\%) & total time (s) &  & OFV & $\mathrm{RPD_b}$ & $\mathrm{RPD_a}$ & best found at (s) \\ \midrule
		10	&&	7405.60	&	7405.60	&	0.00	&	41.29	&&	7405.60	&	0.00	&	0.00	&	0.00	\\
12	&&	8334.20	&	8334.20	&	0.00	&	199.02	&&	8334.20	&	0.00	&	0.00	&	0.00	\\
14	&&	9162.60	&	8966.00	&	1.78	&	1144.58	&&	9154.40	&	-0.02	&	-0.02	&	0.01	\\
20	&&	11828.00	&	8894.00	&	24.56	&	1800.33	&&	11093.80	&	-2.59	&	-2.59	&	0.26	\\
22	&&	12883.80	&	9796.40	&	24.08	&	1800.25	&&	12206.80	&	-4.37	&	-4.37	&	1.27	\\
24	&&	13659.00	&	9207.20	&	32.75	&	1800.12	&&	11850.40	&	-9.78	&	-9.78	&	3.61	\\
50	&&	21040.60	&	12345.20	&	41.27	&	1800.10	&&	18163.00	&	-12.78	&	-11.09	&	18.88	\\
52	&&	21850.20	&	12881.60	&	41.06	&	1800.21	&&	18514.60	&	-14.46	&	-12.52	&	28.95	\\
54	&&	22110.40	&	13029.60	&	40.65	&	1800.12	&&	18786.60	&	-12.98	&	-10.57	&	26.06	\\
80	&&	27023.20	&	16589.60	&	38.72	&	1800.97	&&	24372.80	&	-9.23	&	-6.07	&	57.51	\\
82	&&	25254.00	&	14459.00	&	42.75	&	1802.85	&&	22670.60	&	-10.19	&	-6.73	&	56.43	\\
84	&&	26822.40	&	15222.60	&	43.24	&	1824.55	&&	22898.20	&	-14.37	&	-10.86	&	49.10	\\
100	&&	27199.40	&	15366.40	&	42.54	&	1800.47	&&	23673.80	&	-11.29	&	-8.09	&	67.01	\\
102	&&	29787.20	&	17210.80	&	42.34	&	1801.04	&&	27335.40	&	-8.17	&	-4.56	&	71.90	\\
104	&&	26463.40	&	15653.40	&	40.95	&	1800.12	&&	24744.60	&	-5.76	&	-2.70	&	81.00	\\
\bottomrule
Average	&&	19388.27	&	12357.44	&	30.45	&	1534.40	&&	17413.65	&	-7.73	&	-6.00	&	30.80\\
			\bottomrule
		\end{tabular}%
	}
\end{table}

Results in Tables \ref{tab:Small} -- \ref{tab:Large} show that Gurobi finds optimal solutions within the time limit only when the number of customers is small. In these instances, RKO consistently finds optimal solutions in negligible CPU time. As the number of customers increases, the MILP formulation becomes increasingly complex to solve, as evidenced by the optimality gaps. Across all test sets, the average optimality gap ranges from 22.71\% for instances with $|H|=5$ time intervals to 30.45\% for instances with $|H|=15$ time intervals, even though the solver is provided with an initial incumbent solution. For some combinations of $n$ and $|H|$, the average gap can exceed 42\%.

In contrast, RKO maintains robust performance across all configurations, with more negative relative percentage deviations for both the best and the average solution as the number of time intervals, $|H|$, increases. The average $\mathrm{RPD_b}$ values are equal to $-5.35\%$, $-7.18\%$, and $-7.73\%$ for $|H|=5,10,$ and $15$, respectively, indicating that the best solutions found by RKO are often significantly better than the best upper bounds produced by Gurobi within the time limit. For some combinations of $n$ and $|H|$, the average $\mathrm{RPD_b}$ is greater than $-14\%$. Regarding the average $\mathrm{RPD_a}$, the values are $-3.51\%$, $-5.49\%$, and $-6.00\%$ for $|H|=5,10,$ and $15$, respectively, and for certain instance configurations it is greater than $-12\%$. Moreover, RKO obtains these solutions in a fraction of the CPU time required by Gurobi.

The impact of time discretization is also evident. As $|H|$ increases, the MILP becomes considerably harder to solve, leading to larger optimality gaps and more frequent time-limit violations. In contrast, both the quality of the solutions produced by RKO and the time required to obtain them change only slightly with $|H|$. This shows that the random-key representation combined with the time-dependent decoder is more robust to different time discretizations than the exact formulation, which is more sensitive to the growth in the number of time-indexed decision variables.

Regarding the number of customers, the contrast between the two approaches is clear. For each fixed value of $|H|$, Gurobi reaches the time limit in most instances and stops with substantial remaining gaps, whereas RKO continues to provide solutions with negative $\mathrm{RPD_b}$ and $\mathrm{RPD_a}$ values. Over the 225 TD-TSP instances considered, RKO produced a better upper bound than Gurobi in all but one case: for a single instance with $n = 104$ and $|H| = 5$, the value of $\mathrm{RPD_b}$ was 0.33\%, indicating solutions of very similar quality. Overall, these results suggest that the random-key search combined with the time-dependent decoder is more effective and faster than the exact formulation for solving this problem.

From a modeling and algorithmic standpoint, the results highlight the complementary roles of the MILP formulation and the RKO. The mixed-integer model is useful for obtaining strong bounds and exact solutions on smaller instances, but quickly becomes computationally demanding as the problem size and time discretization grow. RKO, in contrast, encodes only the visiting order of the customers, while the decoder enforces feasibility and evaluates time-dependent travel times over the discretized horizon, allowing the search to focus on the core combinatorial structure. As a consequence, RKO is able to deliver high-quality solutions with predictable computational effort in settings where the exact formulation struggles to produce tight bounds within the time limit.

\begin{figure}[htbp]
    \centering
    \includegraphics[width=0.8\linewidth]{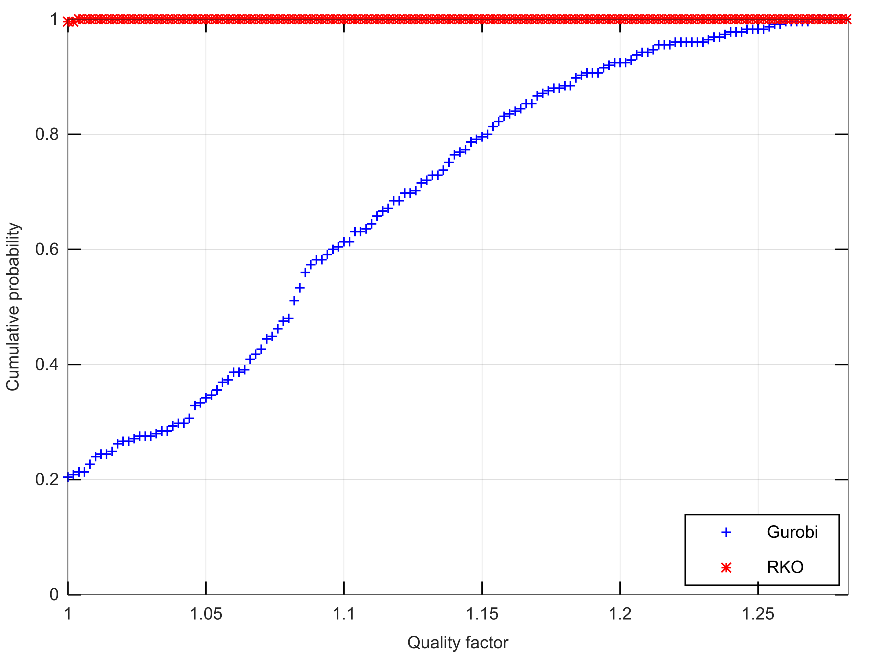}
    \caption{Quality performance profile for the TD-TSP instance set using the best-known value as reference.}
    \label{fig:PP_upper}
\end{figure}

\begin{figure}[htbp]
    \centering
    \includegraphics[width=0.8\linewidth]{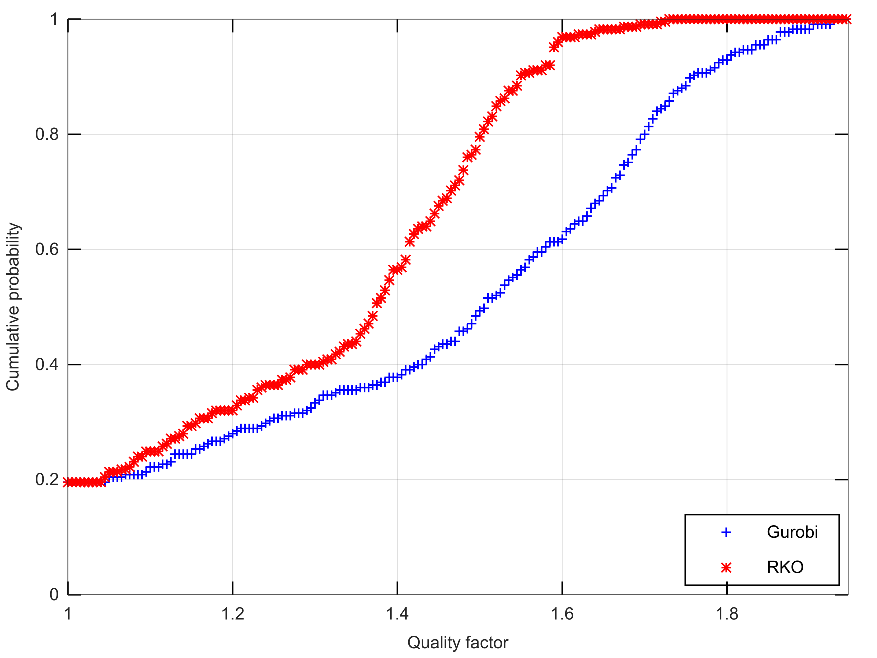}
    \caption{Quality performance profile for the TD-TSP instance set using the Gurobi lower bound as reference.}
    \label{fig:PP_lower}
\end{figure}

To complement the aggregated results above, we finally provide an instance-level comparison using quality performance profiles \cite{dolan2002benchmarking}. These profiles summarize, over all $|P|=225$ instances, the fraction of cases in which each method attains a solution within a given quality factor $\rho$ of a reference value. For each instance $p\in P$, let $z_{p,s}$ denote the best feasible objective value returned by method $s\in\{\text{Gurobi},\text{RKO}\}$ within its prescribed computational budget. Specifically, we take Gurobi’s final incumbent under the 1800s limit and RKO’s best solution across five independent $n$-second runs. We define a quality factor $q_{p,s}\ge 1$ and plot the empirical cumulative distribution $\rho_s(\tau)=\frac{1}{|P|}\,\bigl|\{\,p\in P:\ q_{p,s}\le \tau\,\}\bigr|$ reported as cumulative probability versus the quality factor threshold $\tau$. Two reference schemes are considered. First, we measure quality relative to the best-known solution on each instance, $z_p^{\mathrm{best}}=\min\{z_{p,\text{Gurobi}},z_{p,\text{RKO}}\}$, and $ q_{p,s}^{\mathrm{best}}=\frac{z_{p,s}}{z_p^{\mathrm{best}}}$. Second, we measure quality relative to the Gurobi lower bound $(LB_p)$, using $q_{p,s}^{LB}=1+\frac{z_{p,s}-LB_p}{LB_p}$. Figures \ref{fig:PP_upper} and \ref{fig:PP_lower} report the corresponding profiles.

The performance profiles provide a distributional view that corroborates the table-based analysis. In the profile referenced to the best-known value, RKO attains the best-known value on $224/225$ instances ($\rho_{\text{RKO}}(1)\approx 0.9956$), whereas Gurobi does so on $46/225$ instances ($\rho_{\text{Gurobi}}(1)\approx 0.2044$). Moreover, the markedly slower rise of $\rho_{\text{Gurobi}}(\tau)$ indicates that, under the adopted budgets, its final incumbents are frequently dominated by the best solutions produced by RKO, requiring substantially larger quality ratios to cover the same fraction of instances.

In the profile referenced to the lower-bound value, both methods coincide at $\tau=1$ on $44/225$ instances $\bigl(\rho(1)\approx 0.1956\bigr)$, which corresponds to the subset of cases where the final Gurobi lower bound is tight and the reported primal solution matches the bound. Beyond this subset, RKO remains consistently closer to the bound: for instance, at $\tau=1.5$ it covers $179/225$ instances $\bigl(\rho_{\text{RKO}}(1.5)\approx 0.7956\bigr)$, whereas Gurobi covers only $111/225$ instances $\bigl(\rho_{\text{Gurobi}}(1.5)\approx 0.4933\bigr)$. This dominance indicates that, under the adopted computational budgets, RKO delivers smaller bound-referenced gaps on a larger fraction of instances when both methods are normalized by the same $LB_p$ information. Since the reference values are the final lower bounds produced by Gurobi, this profile should be read as proximity to the best bound that the exact model can certify within the time limit, rather than proximity to the unknown optimum. Moreover, it is not a comparison of lower bounds across methods, but a comparison of the primal upper bounds returned by each method with respect to a common lower-bound reference.

\section{Conclusions}
\label{sec:sec6}

This paper proposes applying the RKO framework to solve MIPs. We evaluated the proposed approach on two combinatorial optimization problems: the canonical mean-variance Markowitz portfolio optimization with buy-in and cardinality constraints, and the Time-Dependent Traveling Salesman Problem (TD-TSP). For each MIP formulation, we develop a problem-specific decoder to handle the constraints and the objective function.

The computational results obtained with RKO are compared with those obtained with the Gurobi solver. In both problems studied, the RKO results are competitive with or superior to Gurobi in terms of solution quality, computational time, and scalability. The RKO outperforms Gurobi on large-scale instances and in 224/225 TD-TSP instances, while remaining competitive on portfolio instances. 

Generally, the random-key vector can represent only the variables with non-zero values. Therefore, RKO solutions reduce the solution space's dimensionality and mitigate flat regions, while enabling the decoder's logic to satisfy as many constraints as possible, thereby generating fewer infeasible solutions that require penalization and facilitating better convergence.


RKO is a viable alternative when commercial solvers become impractical due to computational time or licensing costs. End users can focus entirely on the development of specific decoder designs geared towards other MIP formulations, while the problem-independent search logic is handled out of the box. As demonstrated here, RKO can efficiently handle complex constraints and large instance sizes. However, we stress that the quality of the decoder is critical for obtaining high-quality solutions within a reasonable computational time. 

Future research could explore integrating a parallel Lagrangian relaxation module to complement the RKO framework. Such a hybrid approach would provide dual bounds to assess solution quality and estimate the optimality gap.



We can also extend the application of RKO to other classes of MIPs (e.g., scheduling, lot-sizing, network design). Given its modular nature, RKO can incorporate other metaheuristics or hybrid exact–heuristic strategies, while future research may focus on automated decoder design and theoretical convergence analysis.



\backmatter


\section*{Declarations}

\noindent \textbf{Funding}: This study was funded by FAPESP under grants 2018/15417-8, 2022/05803-3, and 2024/08848-3, and Conselho Nacional de Desenvolvimento Cient\'ifico e Tecnol\'ogico (CNPq) under grants 312747/2021-7 and 405702/2021-3.

\vspace{0.2cm}

\noindent \textbf{Data availability}: The instances and the source code are available at: \url{https://github.com/RKO-solver}

\vspace{0.2cm}

\noindent \textbf{Conflict of Interest}: The authors declare that they have no conflict of interest.

\vspace{0.2cm}


\vspace{0.2cm}

\noindent \textbf{Informed consent} was obtained from all individual participants included in the study.

\bibliography{bibliography}

\end{document}